\documentclass[a4paper,fleqn]{cas-sc}

\usepackage[authoryear]{natbib}
\usepackage{bm}

\usepackage{float}
\usepackage{wrapfig}
\floatstyle{plaintop}
\restylefloat{table}
\usepackage[section]{placeins}
\usepackage[subrefformat=parens,labelformat=parens]{subfig}

\usepackage{tikz}
\usetikzlibrary{shapes.geometric, arrows}
\tikzstyle{startstop} = [rectangle, rounded corners, minimum width=3cm, minimum height=1cm, text centered,draw=black,fill=none]
\tikzstyle{io} = [trapezium, trapezium left angle=70, trapezium right angle=110, minimum width=1cm, minimum height=1cm, text centered, draw=black, fill=none]
\tikzstyle{process} = [rectangle, minimum height=1cm, text width=6cm, text centered, draw=black, fill=none]
\tikzstyle{decision} = [diamond, minimum width=3cm, minimum height=1cm, text centered, draw=black, fill=none]
\tikzstyle{arrow} = [thick,->,>=stealth]

\usepackage{natbib}
 \bibpunct[, ]{(}{)}{,}{a}{}{,}%
 \def\bibsep{\smallskipamount}%
\usepackage{placeins}
\usepackage{float}
\usepackage{bm}
\usepackage{wrapfig}
\usepackage{subfig}
\usepackage{caption}
\usepackage{multirow}
\usepackage{vcell}
\usepackage{booktabs}
\usepackage[export]{adjustbox}
\usepackage{graphicx, booktabs, multirow, makecell}
\usepackage{tabularx}
\newcolumntype{Y}{>{\centering\arraybackslash}X}

\usepackage{amsmath}
\usepackage{amsthm}
\usepackage{subfig}
\usepackage{graphicx}
\usepackage{setspace}
\usepackage{multicol}
\usepackage[ruled,vlined]{algorithm2e}
\SetKwInput{KwInput}{Input}                % Set the Input
\SetKwInput{KwOutput}{Output}              % set the Output
\SetKwProg{try}{try}{:}{}
\SetKwProg{catch}{catch}{:}{end}

\newcolumntype{M}{>{\centering\arraybackslash}m{2cm}}
\newcolumntype{N}{>{\centering\arraybackslash}m{0.6cm}}
\usepackage{xcolor}

\def\tsc#1{\csdef{#1}{\textsc{\lowercase{#1}}\xspace}}
\tsc{WGM}
\tsc{QE}
\tsc{EP}
\tsc{PMS}
\tsc{BEC}
\tsc{DE}
\begin{document}
\def\floatpagepagefraction{1}
\def\textpagefraction{.001}
\shorttitle{T. Lilasathapornkit et~al.}
\shortauthors{T. Lilasathapornkit et~al.}
%\begin{frontmatter}

\title [mode = title]{Dynamic Pedestrian Traffic Assignment with Link Transmission Model for Bidirectional Sidewalk Networks}

\author[1]{Tanapon Lilasathapornkit}

%\ead{t.lilasathapornkit@unsw.edu.au}

\address[1]{Research Centre for Integrated Transport Innovation (rCITI), School of Civil and Environmental Engineering, University of New South Wales, Sydney, NSW, Australia}

\author[1]{Meead Saberi}
\cormark[1]

%\ead{meead.saberi@unsw.edu.au}

\begin{abstract}
Planning assessment of the urban walking infrastructure requires appropriate methodologies that can capture the time-dependent and unique microscopic characteristics of bidirectional pedestrian flow. In this paper, we develop a simulation-based dynamic pedestrian traffic assignment (DPTA) model specifically formulated for walking networks (e.g. sidewalks) with bidirectional links. The model consists of a dynamic user equilibrium (DUE) based route choice and a link transmission model (LTM) for network loading. The formulated DUE adopts a pedestrian volume delay function (pVDF) taking into account the properties of bidirectional pedestrian streams such as self-organization. The adopted LTM uses a {\color{black}three-dimensional} triangular bidirectional fundamental diagram as well as a generalized first-order node model. The applicability {\color{black}and validity} of the model is demonstrated in hypothetical small networks as well as a real-world large-scale network of sidewalks in Sydney. {\color{black}{}The model successfully replicates formation and propagation of shockwaves in walking corridors and networks due to bidirectional effects.}

\end{abstract}


\maketitle

\section{Introduction}

Potential health and societal benefits of active transportation are becoming more acknowledged. Many cities around the world are increasing their investment in walking infrastructure. However, overcrowded footpaths in some cities during peak hours create potential safety risks and increase delays for pedestrians. Investment in walking infrastructure is often made on an adhoc manner rarely supported by strategic large-scale pedestrian network models similar to what is commonly used for analysis of vehicular traffic systems.

Urban footpaths or sidewalks can be viewed as a network of bidirectional pedestrian links. Several studies in the past have already investigated the bidirectional crowd dynamics using the fundamental relationship between flow and density \citep{Seyfried2005,zhang2012ordering,hanseler2014macroscopic,cao2017fundamental,hanseler2017dynamic, saberi2014, saberi2015}. Despite its significance and practical relevance, very little effort has been put into understanding and modeling the network-wide impact of pedestrian traffic in the urban context for planning applications. This study aims to develop a simulation-based dynamic pedestrian traffic assignment (DPTA) framework to model large-scale footpath or sidewalk networks.

Research on pedestrian flow modeling and dynamics has grown in several directions in the past few decades including development of novel approaches in microscopic modeling \citep{lovaas1994modeling,helbing1995social,black2001cellular,moussaid2010walking,huang2017behavior,tao2017cellular,shahhoseini2018pedestrian}, mesoscopic modeling \citep{xiong2010hybrid,cristiani2011multiscale,tordeux2018mesoscopic}, and macroscopic simulation \citep{hughes2002continuum,colombo2011non,schwandt2013macroscopic,hoogendoorn2015continuum,hanseler2017dynamic,hoogendoorn2018macroscopic,taherifar2019macroscopic,aghamohammadi2020dynamic,moustaid2021macroscopic, molyneaux2021design}. For a detailed review of the literature on pedestrian crowd data collection and empirical insights, see \citep{haghani2018crowd,haghani2020empirical1,haghani2020empirical2}. 

Simulation-based dynamic traffic assignment (DTA) models are used as tools to provide estimates of traffic conditions for transportation operations applications and are increasingly applied in strategic transportation planning to support decision making of infrastructure investments. Two main reasons that have led to increasing popularity of DTA developments in practice in the past few decades are as follows. Firstly, the static traffic assignment model does not sufficiently capture congestion dynamics in the network including oversaturation, formation of physical queues, spillbacks, and shockwaves \citep{raadsen2016efficient, bliemer2019continuous}. In a previous study by \cite{lilasathapornkit2020traffic}, a static traffic assignment modeling framework with bidirectional link performance functions is already proposed. Secondly, implementation of a sub-daily traffic model requires anticipation of spatial and temporal changes in traffic conditions that could result in re-routing or traffic accumulation in certain areas of the network \citep{melson2018dynamic}. Simulation-based DTA models capture the essential dynamic properties of traffic flow including planned interruptions such as signal control \citep{zhu2015linear} or unplanned disruptions such as capacity reduction due to incidents \citep{he2012modeling,nogal2016resilience}. A large body of research already exists that looks into increasing the efficiency and accuracy of traffic simulation models in large-scale networks while retaining the most essential characteristics of traffic flow \citep{saeedmanesh2017dynamic,yildirimoglu2018hierarchical,batista2019regional,han2019computing}. \citet{aghamohammadi2020dynamic} provides a comprehensive literature review on DTA modeling for both pedestrians and vehicles. In this study, we build upon the rich literature of simulation-based DTA and develop the first dynamic sidewalk network model that takes into account the microscopic behaviour of bidirectional pedestrian traffic flow such as self-organization.

 Network traffic user equilibrium (UE) problem can be formulated as Variational Inequality (VI) in the context of static traffic assignment \citep{dafermos1980traffic} and dynamic networks \citep{friesz1993variational, ran1996link, bliemer2001analytical, friesz2011approximate, gentile2016solving}. However, the main difficulty with the analytical DTA approach is adding realistic traffic dynamics to an already sophisticated formulation. Reproducing the formation, propagation , and dissipation of physical queues across multiple links are not a trivial task. The seminal LWR model uses the kinematic wave theory to address this limitation \citep{lighthill1955kinematic, richards1956shock}. To solve the LWR problem, \citet{Daganzo1994,daganzo1995cell} proposed the cell transmission model (CTM) that discreticizes time and space and provides approximate solutions. A few years later \citet{yperman2005link} proposed the link transmission model (LTM) that discriticizes time and link boundary conditions based on Newell's simplified formulation of kinematic wave theory \citep{newell1993simplified}. LTM often yields a more accurate solution than CTM with the same time-step size and less computational cost \citep{jin2015continuous}. During the past decade, LTM has become an increasingly popular model for dynamic network loading (DNL) \citep{gentile2010general, himpe2016efficient, raadsen2016efficient, Chakraborty2018, bliemer2019continuous, raadsen2019continuous}. In this study, we adopt a similar LTM approach and integrate it into a simulation-based DTA to model footpath networks with bidirectional links.

{\color{black}The main contributions of this study are as follows: (i) We formulate a DUE problem that adopts a travel cost function appropriate for bidirectional pedestrian streams calibrated with data from controlled experiments; (ii) We propose and calibrate new bidirectional three-dimensional pedestrian fundamental diagrams for traffic assignment applications; (iii) We also propose a new node model to accommodate for the bidirectional pedestrian flows within the LTM; (iv) We demonstrate the applicability and validity of the model in both small and large-scale pedestrian networks. The paper provides a balance between mathematical rigor on one hand and realism on the other. The proposed pedestrian simulation-based DTA framework adopts the classical DUE with flow-based cost functions for the route choice model and the LTM for DNL. To the best knowledge of the authors, this is the first study in the literature that presents a large-scale LTM-based DPTA model for bidirectional sidewalk networks.}

The remainder of the paper is organized as follows. A review of the literature on route choice models and DNL is provided in Section \ref{sec:routechoice} and Section \ref{sec:dnl}. Section \ref{sec:example} presents the numerical results including a small grid network, a long corridor and a large-scale sidewalk network of Sydney. Section \ref{sec:conclusion} summarizes the contributions and suggests future research directions.

\section{Methodology}
\begin{figure}
	\centering
	\includegraphics[width=0.9\linewidth]{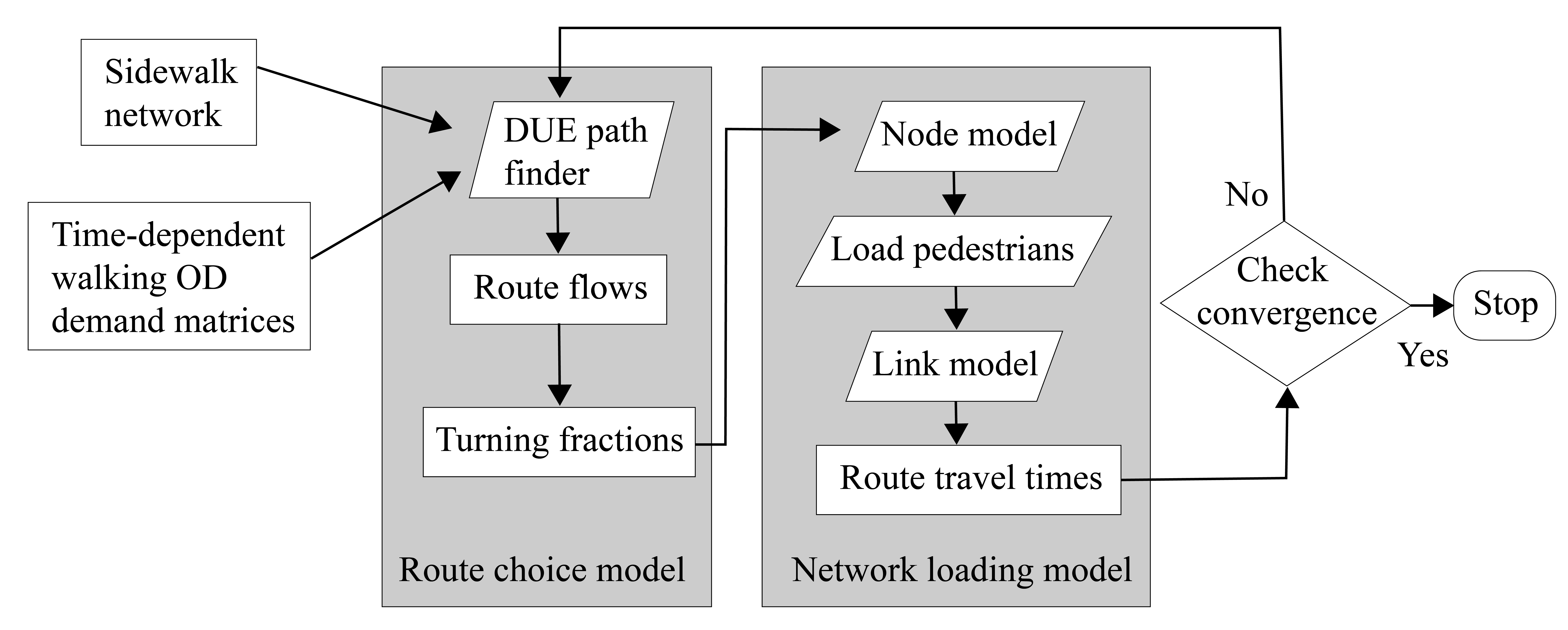}
	\caption{An overview of the proposed DPTA model.} 
	\label{fig:dtaoverview}
\end{figure}

A DTA model consists of two general components: the
route choice model and the DNL. In the following, we describe the methodological details of both components in our developed approach as illustrated in Figure \ref{fig:dtaoverview}. 

\subsection{Route choice model} \label{sec:routechoice}
Link travel times vary over time depending on the prevailing traffic conditions. Link travel times can be estimated based on link inflow \citep{jang2005discrete,bliemer2003quasi,gentile2016solving}, occupancy \citep{lo2002cell,daganzo1995properties}, or horizontal difference between cumulative counts \citep{tong2000predictive,nie2010solving,friesz2013dynamic}. In this study, we will use an inflow based link travel time function extended from our previous work on pedestrian volume delay functions \citep{lilasathapornkit2020traffic}.

In the early DTA models, equilibrium conditions were often based on instantaneous path cost \citep{friesz1993variational} as in static traffic assignment models. However, more recent implementations of DTA often include experienced path cost instead \citep{chiu2011dynamic}. In this study, travellers are assumed to have a perfect knowledge of the network including current travel costs (instantaneous path cost), but do not have the knowledge of their actual travel cost (experienced path cost). Assignment of traffic based on instantaneous route travel times is not necessarily incorrect, but it is under assumption that travellers decide to take the shortest route using pre-trip travel information \citep{ran1993new,chiu2011dynamic,chen2014real, yildirimoglu2013experienced}. 

A DUE algorithm iteratively updates route choices over the entire network until all users going from origin $r$ to destination $s$ at time $t$ share the same minimal expected travel time. {\color{black}Solving DUE with continuous time variable may not be practical in large-scale networks, so discrete time formulation is often utilized \citep{chen1998model,bliemer2003quasi}}. \citet{yperman2005link} introduced the first LTM based on a simplified kinematic wave theory with a discrete time. Several studies in the past aimed to reduce LTM computation time by proposing efficient iterative algorithms \citep{himpe2016efficient}, event-based algorithms \citep{raadsen2016efficient}, and continuous algorithms \citep{han2016continuous,gentile2016solving,raadsen2019continuous}. In this study, we use a discrete-time approach.

Building upon the Wardop's first principle \citep{Wardrop1952}, we use the DUE condition in which travelers route choice is based on the minimum instantaneous route costs. However, the DNL, as will be described later in the paper, estimates link and route travel times based on a triangular fundamental diagram that estimates the actual (experienced) route travel costs. The two route travel time estimates are expected to be highly correlated as shown later in the numerical experiments section of the paper.

\subsubsection{Dynamic User Equilibrium.}

{\color{black}
Let $G=(N,A)$ be the directed graph consisting of nodes $N$ and links $A$. Let \textbf{c} be the vector of link travel time functions, \textbf{u} be the vector of link flows, $\bar{u}$ be the vector of optimized link flow, $S$ be a set of feasible link flow patterns, $K$ be a set of discrete departing time instants, $T$ be a set of discrete time instants, $u_a$ be the pedestrian flow on link $a \in A$, $N$ be a set of nodes, $\Pi_{rs}$ denotes the set of paths connecting OD pair $(r,s) \in W$. Let $c_{a,t} (u_{a,t})$ be the walking travel time on link $a \in A$ as a function of the link flow $u_{a,t}$ at time instance{t}, $c_{p,k} ^{rs}$ be the path travel time on path $p$ at the departure time $k \in K$, $\delta^{rs}_{a,p,k,t}$ be the link-path incidence matrix at time $t$, $f^{rs}_{k}$ be the flow on path $p$ at the departure time $k$ connecting an OD pair $(r,s) \in W$, and $q_{k}^{rs}$ be the travel demand between OD pair $(r,s)$ at the departure time $k$.  \citet{Wardrop1952}'s first principle states that at UE, no traveler has any incentive to unilaterally change its route. Therefore, the traffic assignment problem (TAP) can be formulated as a variational inequality problem VI(\textbf{c},S) that finds the optimal link flows $\bar{u}$ as expressed in Equation~\eqref{eq:wardrop1st}. 
}
\begin{subequations}\label{eq:wardrop1st} \color{black}
\begin{align}
& && \textbf{c}  (\bar{u} )(\textbf{u}  - \bar{u} ) \geq 0 && \forall \textbf{u} , \bar{u} \in S \label{eq:wardrop1st_a} \\
&\textrm{subject to:} && u_{a,t} = \sum_{r \in N} \sum_{s \in N} \sum_{p \in \Pi_{rs}} \sum_{k \in K} \delta^{rs}_{a,p,k,t} f^{rs}_{p,k}  &&\forall a \in A, \forall t \in T \label{eq:wardrop1st_b}\\
&&&\sum_{p \in \Pi_{rs}} f^{rs}_{p,k} = q_{k} ^{rs} &&\forall (r,s) \in W, \forall k \in K \label{eq:wardrop1st_c}\\
&&& f^{rs}_{p,k} \geq 0 &&\forall (r,s) \in W, \forall p \in \Pi_{rs}, \forall k \in K \label{eq:wardrop1st_d}\\
&&& c_{p,k}^{rs} = \sum_{a \in A} \sum_{t \in T} c_{a,t} \delta^{rs}_{a,p,k,t} && \forall p \in \Pi_{rs}, \forall k \in K \label{eq:wardrop1st_e}
\end{align}
\end{subequations}

\subsubsection{Link travel time function}
{\color{black}Here, we use our previously developed pedestrian volume delay function (pVDF) in \citet{lilasathapornkit2020traffic}.}
Consider a pair of link $a, a' \in A$ that belongs to the same bidirectional stream. Let $\tau _{a}$ denotes the free flow travel time of link $a$, {\color{black}$C_{a}$ denotes the capacity of link $a$}, and $\alpha$ and $\beta$ be model parameters. We denote 	
\begin{equation} \color{black}
c_{a,t}  \left(u_{a,t} ,u_{a',t}\right) = \tau _{a} \left( 1+\alpha \left(\frac{u_{a,t}  +u_{a',t} }{C_{a}}\right)^{\beta} \right) \hspace{120pt} \quad \forall a, a' \in A, \forall t \in T
\label{eq:detvdf_sym}
\end{equation}    	
{\color{black}the deterministic symmetric pVDF of link $a \in A$ for time $t \in T$.}
Equation \eqref{eq:detvdf_sym} represents a pVDF with symmetric link interactions that is extended from the well-known BPR function with $\alpha$ representing the ratio of travel time per unit distance at practical capacity to that at free flow while $\beta$ determines how fast the average travel time increases from free-flow to congested conditions. Links from both directions of the same bidirectional stream have identical travel times. Equation \eqref{eq:detvdf_sym} also implies that flows from either directions have equivalent impact on travel time of both directions. 

We also consider an asymmetric pVDF as follows where $\alpha, \beta, \mu, \eta _{r}, \lambda _{r}, \eta _{c},$ and $ \lambda _{c}$ are the model parameters. We denote

\begin{equation} \label{eq:detvdf_asym} \color{black}
 c_{a,t}  (u_{a,t} ,u_{a',t} ) = \tau _{a} \Big ( 1+\alpha \Big(\frac{u_{a,t}+u_{a',t} }{C_{a}}\Big)^{\beta} \Big ) + \tau _{a} \mu \mathrm{e} ^{\eta _{r} \big ( \frac{u_{a,t} }{C_{a}}- \lambda _{r} \big ) ^{2} + \eta _{c} \big ( \frac{u_{a',t} }{C_{a}}- \lambda _{c} \big ) ^{2}} \hspace{10pt} \quad \forall a, a' \in A, \forall t \in T
\end{equation}

{\color{black}the deterministic asymmetric pVDF of link $a \in A$ for time $t \in T$.}

The first term in Equation \eqref{eq:detvdf_asym} is the same as that in Equation \eqref{eq:detvdf_sym}. The second term, referred to as the bidirectional term, captures the bidirectional pedestrian flow in which as pedestrian flows get closer to a certain value ($\eta _{r}$ for the reference direction and $\eta _{c}$ for the counter direction), the travel time substantially increases. The bidirectional term becomes negligible when flows are further away from these values ($\eta _{r}$ and $\eta _{c}$). The term consists of five parameters. $\mu$ determines the magnitude of the bidirectional impact for the balanced flows. $\eta _{r}$ and $\eta _{c}$ determine the range of flows that exhibit bidirectional impact by regulating the width of a bell-shaped curve base. $\lambda _{r}$ and $\lambda _{c}$ determine the flow ratios that have the highest congestion level in the stream. {\color{black}In this study, we use the proposed link travel time functions in \cite{lilasathapornkit2020traffic} as the cost functions in the route choice model.} The estimated travel time is treated as a pre-trip expected travel time used to find the DUE. 

\subsubsection{Turning fraction}
The distribution of flows from incoming links toward outgoing links is obtained from the total turning fraction. Let $\phi _{ij} ^{n}$ be the turning fraction from incoming link $i$ toward outgoing link $j$ through node $n \in N$. We denote
\begin{subequations}\label{eq:turningfraction}
\begin{align}
 & 0 \leq \phi _{ij} ^{n} \leq 1 , \hspace{50pt} && \forall i,j \in A , \forall n \in N \\
& \sum _{j} \phi _{ij} ^{n} = 1, && \forall i \in A , \forall n \in N
\end{align}
\end{subequations}
We derive the turning fractions from the output of the DUE by converting all route flows toward the same destination into a set of turning fractions at each time instant. For example if we have 3 origins and 2 destinations, 2 sets of turning fractions will be determined.

\subsection{Dynamic network loading} \label{sec:dnl}
A first order DNL model consists of a link model and a node model often based on kinematic wave theory  \citep{lighthill1955kinematic,richards1956shock} in which vehicles are represented as flows that propagate in the network. The link model defines flow propagation within each link boundary while the node model regulates how upstream links distribute flows toward downstream links at each node. Here, we propose a pedestrian triangular fundamental diagram that captures bidirectional interactions in walking streams. The DNL model loads users into associated routes though the network from the origin nodes towards the destination nodes. We use LTM as our DNL model due to its calculation efficiency \citep{yperman2005link}. LTM consists of the flow model which follows a modified kinematic wave theory \citep{newell1993simplified} and the node model that utilizes route choice to depict flow propagation throughout the network. Network users may encounter a different experienced travel time than the expected travel time as a result of the DNL. If more than one destination are present, a multi-commodity DNL is utilized.

Two main causes could create a pedestrian traffic shockwave. First, a reduction in the link width leads to reduction in link capacity and as a result, a shockwave forms. In freeway traffic, a sudden capacity drop may lead to shockwave formation \citep{wu2011shockwave}. Similarly a shockwave can emerge in unidirectional pedestrian traffic as well \citep{helbing1995social, zhang2013empirical}.  In this study, shockwave speed is constant and does not vary based on the free flow speed. 

Secondly, pedestrian flow from one direction decreases the effective capacity of the other direction. This phenomena isn't yet fully understood and has only been explored in a few past studies \citep{zhang2015simulation,feliciani2018universal,fujita2019traffic}.

\subsubsection{Link Model}

\citet{long2011discretised} proposed link travel time functions that are flow dependent as either a step function or a linear interpolation function. The function, however, does not capture the spillback effect. An alternative is using the fundamental relationship between flow and density. \citet{flotterod2015bidirectional} proposed a pedestrian bidirecional fundamental diagram that was calibrated with empirical data. The fundamental diagram proposed by \citet{flotterod2015bidirectional} provides a good representation of pedestrian bidirectional dynamics. However, it cannot be used in the original LTM framework \citep{yperman2005link,yperman2007kinematic} as the original LTM requires a triangular fundamental diagram \citep{newell1993simplified} to keep track of the traffic state. Several later studies extended the LTM formulations with piece-wise linear and concave FDs \citep{yperman2007kinematic,gentile2010general, van2017extending,bliemer2019continuous, raadsen2019continuous}.

\begin{figure}
	\centering
	\includegraphics[width=0.8\linewidth]{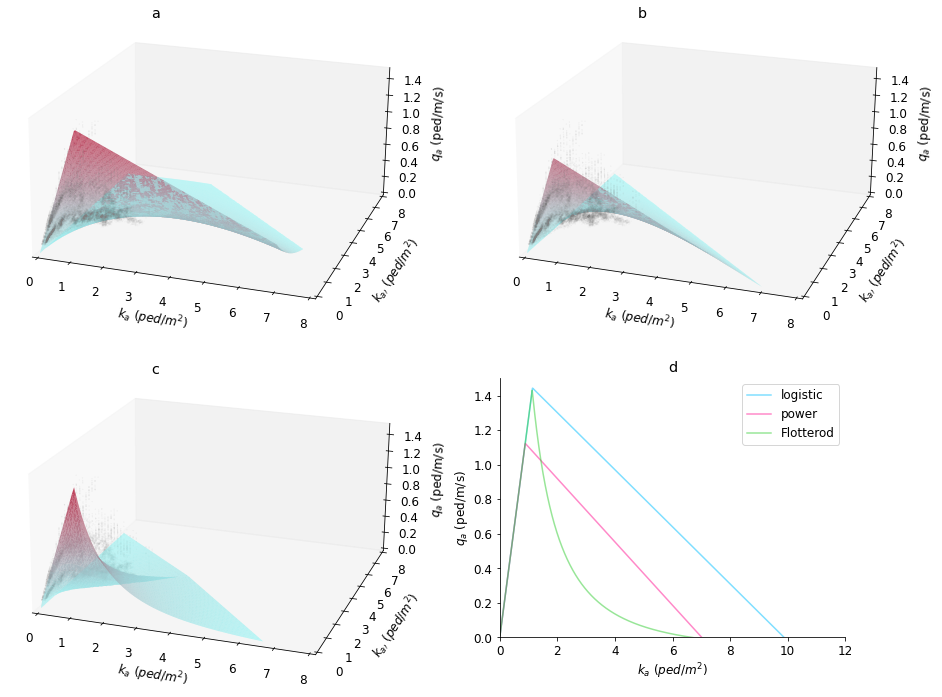}
	\caption{\color{black}Three dimensional bidirectional fundamental diagram based on (a) a logistic function (b) a power function (c) the function proposed by \citet{flotterod2015bidirectional}. Grey points represent empirical data used for calibration (d) A comparison between three FDs for unidirectional flow when density ratio $\rho _{a}$ = 1.} 
	\label{fig:triangularFD_multiple}
\end{figure}

Here, we use a triangular FD and proposes a three-dimensional extension to bidirectional pedestrian streams similar to what was previously proposed in \citet{flotterod2015bidirectional}. The proposed triangular bidirectional FD can be constructed from free flow speed $v_{f}$, shock wave speed $\omega$, and jam density $k_{j}$ as expressed in Equation \eqref{eq:density_ratio} - \eqref{eq:qk_triangular}. {\color{black} The density ratio $\rho _{a}$ represents the ratio of the density of the reference direction over the density of the opposite direction in a bidirectional stream. The effective jam density $\hat{k}_{jam,a}$ represents the practical jam density of the reference direction considering the density from the opposite direction. Similarly, the effective free flow speed $\hat{v}_{f,a}$ represents the free flow speed considering the bidirectionality effects.} Consider a pair of links $a, a' \in A$ which belong to the same bidirectional stream. Let $k _{a}$ and $k _{a'}$ be density of link $a$ and $a'$, respectively. Let $\rho _{a}$ be the density ratio expressed as
\begin{equation} \label{eq:density_ratio}
 \rho_{a} = \frac{k_{a}}{k_{a}+k_{a'}}
\end{equation}
Let $\hat{k}_{jam,a}$ be the effective jam density which is proportional to density ratio and can be defined as
\begin{equation} \label{eq:density_jam}
 \hat{k}_{j,a} = \rho_{a} k_{j,a} 
\end{equation}
Let $k_{c,a}$ be the critical density defined as
\begin{equation} \label{eq:density_critital}
 k_{c,a} = \frac{\hat{k}_{j,a} \omega _{a} }{ \hat{v}_{f,a} + \omega _{a}}
\end{equation}
where $\omega _{a}$ denotes the shockwave speed of link $a$ and $\hat{v}_{f,a}$ denotes the effective free flow speed of link $a$ that can be defined in different ways as
\begin{subequations}\label{eq:freeflowspeed}
\begin{align}
 & \hat{v}_{f,a} = \frac{ v_{f,a} }{{\rm e}^{1 - \rho_{a}} } \label{subeq:freeflowspeed_logistic} \\
& \hat{v}_{f,a} = \rho_{a}^{\gamma}  v_{f,a} \label{subeq:freeflowspeed_power} 
\end{align}
\end{subequations}
where $v_{f,a}$ denotes the free flow speed of link $a$, $\gamma$ denotes the model parameter for calibration. Equation \eqref{subeq:freeflowspeed_logistic} expresses the relationship between density ratio and free flow speed as a logistic function while Equation \eqref{subeq:freeflowspeed_power} expresses the relationship as a power function. The density ratio decreases if the number of pedestrians from the direction of interest decreases or the number of pedestrians from the opposite direction increases. As density ratio reduces, effective speed becomes slower which is reflected in Equation \eqref{eq:freeflowspeed}. However, for unidirectional flows, the effective speed would be identical with free flow speed and both Equation \eqref{subeq:freeflowspeed_logistic} and Equation \eqref{subeq:freeflowspeed_power} would be equivalent. The reduction in speed suggests that pedestrians experience a lower manoeuvrability and freedom to move as higher number of pedestrians flow from the opposite direction. In this study, we assume that the shockwave speed remains constant independent of the density ratio. Figure \ref{fig:triangularFD_multiple} provides a comparison between the FD proposed by \citet{flotterod2015bidirectional} and the proposed three dimensional triangular bidirectional FD in this study for use in the LTM. {\color{black} The three-dimensional triangular bidirectional FDs are calibrated using pedestrian trajectory data from a set of controlled bidirectional experiments \citep{zhang2012ordering, boltes2013} with up to 350 participants. Pedestrian traffic density and flow measurements are then estimated using extended Edie’s definitions to three-dimensional time-space diagrams \citep{saberi2014}.  

An upstream link boundary can either be at free flow state or congested state, implying either hypocritical or hypercritical flow state in the fundamental diagram. Based on the traffic flow fundamental identity $q = v k$, we express flow and density relationship as follows
\begin{subequations}\label{eq:qk_triangular}
\begin{align}
 & q_{a} (k_{a}) = \hat{v}_{f,a} k_{a}  \hspace{50pt}&& \text{if } k_{a} \leq k_{c,a} \label{subeq:qk_triangular_hypo} \\
& q_{a} (k_{a}) = \omega _{a} ( \hat{k}_{jam,a} - k_{a} ) && \text{if } k_{a} > k_{c,a} \label{subeq:qk_triangular_hyper} 
\end{align}
\end{subequations}
The logistic function has its capacity closer to the capacity in \citet{flotterod2015bidirectional} FD, while the power function has a closer jam density. It is more important to have a closer capacity than jam density since a shockwave occurrence is less common than a capacitated flow. The jam density is utilized less frequently than sending flow is at its capacity. So, we decided to use the logistic function in the case studies.}

\FloatBarrier

{\color{black}\citet{yperman2005link} proposed a numerical solution for the LTM that estimates the sending flow $S_{i}$ of an incoming link $i$ and the receiving flow $R_{j}$ of an outgoing link $j$. Equation \eqref{subeq:sendingflow_a} expresses the upstream boundary condition at time $t$ as the outflow of incoming link $i$ to a node between time $t$ and $t + \Delta t$. Link capacity condition is defined as the maximum flow that can be sent from incoming link $i$ limited to its capacity as expressed in Equation \eqref{subeq:sendingflow_b}. Sending flow $S_{i}(t)$ is restricted by both conditions as expressed in Equation \eqref{subeq:sendingflow_c}.}

\begin{subequations}\label{eq:sendingflow}
\begin{align}
 & S_{i,boundary} (t) = U_{i}(t + \Delta t - L_{i}/ \hat{v}_{f,i}) - V_{i}(t) \label{subeq:sendingflow_a} \\
& S_{i,link}(t) = C_{i} \Delta t \label{subeq:sendingflow_b} \\
& S_{i} (t) = min(S_{i,boundary}(t), S_{i,link}(t)) \label{subeq:sendingflow_c}
\end{align}
\end{subequations}
{\color{black} where $U_{i}(t)$ denotes the cumulative number of pedestrians that pass through upstream of incoming link $i$ at time $t$.}
Similarly the downstream boundary condition at time $t$ is the inflow of outgoing link $j$ from a node between time $t$ and $t + \Delta t$ as expressed in Equation \eqref{subeq:receivingflow_a}. Link capacity condition is defined as the maximum flow that can be sent from outgoing link $j$ limited to its capacity as expressed in Equation \eqref{subeq:receivingflow_b}. Receiving flow $R_{j}(t)$  is restricted by both conditions as expressed in Equation \eqref{subeq:receivingflow_c}.

\begin{subequations}\label{eq:receivingflow}
\begin{align}
 & R_{j,boundary} (t) =  V_{j}(t + \Delta t - L_{j}/ \hat{v}_{w,j}) + k_{j,j} L_{j} - U_{j}(t) \label{subeq:receivingflow_a} \\
& R_{j,link}(t) = C_{j} \Delta t \label{subeq:receivingflow_b} \\
& R_{j} (t) = min(R_{j,boundary}(t), R_{j,link}(t)) \label{subeq:receivingflow_c}
\end{align}
\end{subequations}
{\color{black} where $V_{j}(t)$ denotes the cumulative number of pedestrians that pass through downstream of outgoing link $j$ at time $t$.}

Both sending and receiving flows are constraints to limit how many flows can travel from/to any link. In the next subsection, actual flows are regulated by a node model to maximize the total flow $q_{ij} ^{n}$ from incoming links $i$ to outgoing links $j$ through node $n \in N$  while satisfying the constraints.

\subsubsection{Node Model}
Several studies in the literature have explored node merging and diverging behaviors \citep{daganzo1995cell,bliemer2007dynamic,gentile2010general,smits2015family}. Here, we follow a generic class of first order macroscopic node model formulation \citep{tampere2011generic} that transfers total flows consistently over all nodes in the network \citep{himpe2016efficient,raadsen2016efficient}. Let $N$ be a set of nodes in the network and $q_{ij}^{n}$ be the flow from incoming link $i$ toward outgoing link $j$ through node $n$. Let $S _{i} ^{n}$ be the demand constraint as the maximum flow that incoming link $i$ could possibly send if the node $n$ and outgoing link(s) impose no restriction. Let $R _{j}$ be the supply constraint as the maximum flow that outgoing link $j$ could {\color{black}possibly} receive if the node $n$ and incoming link(s) impose no restriction. Let $\phi _{ij} ^{n}$ be the turning fraction of node $n$ from the incoming link $i$ toward the outgoing link $j$ which is defined as the ratio of the total flow from the incoming link $i$ to the outgoing link $j$ over total incoming flows from link $i$. $j'$ is the incoming link of the same bidirectional stream. The node model is then formulated as a maximization problem as follows:

\begin{subequations}\label{eq:nodemodel}
\begin{align}
   & \textrm{maximize} &&  \sum _{i} \sum _{j} q_{ij} ^{n} &&&\forall n \in N \label{subeq:nodemodel_a} \\
   & \textrm{subject to: }  && q_{ij} ^{n} \geq 0 , &&&\forall i,j \label{subeq:nodemodel_b}\\
   & && q_{i} ^{n} = \sum _{j} q_{ij} ^{n} \leq S_{i} ^{n}, &&&\forall i \label{subeq:nodemodel_c}\\
   & \hspace{40pt} q_{j} ^{n} + \tilde{S}_{j} ^{n} = &&\sum _{i} q_{ij} ^{n} + \Big ( U_{j'}( t + \Delta t - L_{j}/v_{f,i}) - U_{j'}( t - L_{j}/v_{f,i} ) \Big ) \leq R_{j} ^{n}, &&&\forall j, j' \label{subeq:nodemodel_d}\\
   &  && \phi _{ij} ^{n} = \frac{S_{ij} ^{n}}{S_{i} ^{n}} = \frac{q_{ij} ^{n}}{q_{i} ^{n}}, &&&\forall i,j \label{subeq:nodemodel_e}
\end{align}
\end{subequations}

Equation \eqref{subeq:nodemodel_b} imposes a restriction that traffic would never flow backward as a non-negativity constraint. Conservation of flow is ensured with Equation \eqref{subeq:nodemodel_c} and \eqref{subeq:nodemodel_d} as the demand and supply constraints, respectively. The supply constraint consists of two terms, $\sum _{i} q_{ij} ^{n}$ represents the sum of all incoming link flows going toward outgoing link $j$. While the other term represents inflow of incoming links $j'$ which is on the same bidirectional stream as outgoing link $j$, taking into account the \lq\lq look ahead" behavior of the opposite direction prior to reaching node $n$. Equation \eqref{subeq:nodemodel_e} enforces conservation of turning fractions (CTF) which ensures the ratio of partial demand $S_{ij} ^{n}$ over the total demand $S_{i} ^{n}$ is equal to the ratio of partial flow $q_{ij} ^{n}$ over the total flow $q_{i} ^{n}$. CTF condition is equivalent to the first-in-first-out (FIFO) condition \citep{tampere2011generic}. If the flow is restricted by one of the outgoing links, the total flows will be the scaled-down version of the total demand. This condition prevents node model from unrealistically preferring any flows over others. Furthermore, the invariance principle condition is automatically satisfied since $q_{ij}$ is derived by distributing supply and not the demand \citep{tampere2011generic}.

{\color{black}Overall, the node model regulates how many pedestrians pass through node $n$ at time $k$ after pedestrians already decide which paths to take. However, pedestrians that are traveling in the network may adjust their path based on dynamically changing traffic conditions.}

\begin{figure}
	\centering
	\includegraphics[width=0.2\linewidth]{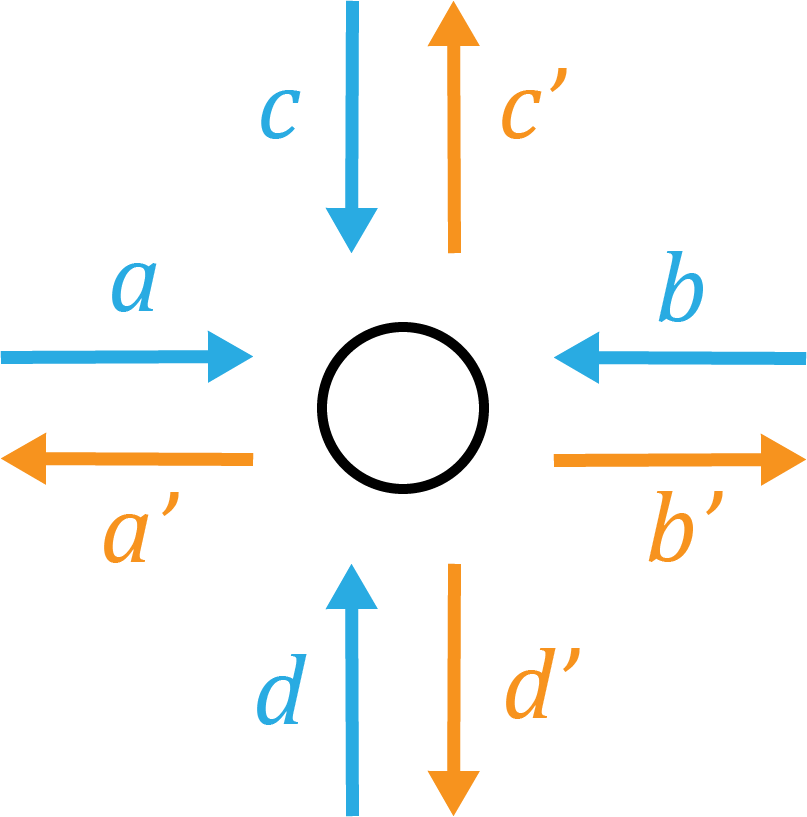}
	\caption{A node model example} 
	\label{fig:nodemodel_example}
\end{figure}

Here, we provide a simple example to demonstrate how the node model can be applied. A standard intersection consisting of four incoming and four outgoing links is considered as shown in Figure \ref{fig:nodemodel_example}. From the given demand in Table \ref{tab:nodemod1}, the total flows $q_{ij} ^{n}$ are restricted by several constraints as expressed in Equations \eqref{subeq:nodemodel_b} - \eqref{subeq:nodemodel_e}.

\FloatBarrier

\begin{table}
\centering
\caption{Demand $S_{ij}$}
\label{tab:nodemod1}
\centering
    \begin{tabularx}{\textwidth}{Y | YYYY | Y}
    \hline
    $S_{ij}$   & $a'$  & $b'$  & $c'$  & $d'$ & $S_{i}$   \\ \hline
    $a$   & 0   & 1   & 0   & 0    & 1   \\
    $b$   & 1   & 0   & 0   & 0.5  & 1.5 \\
    $c$   & 0   & 0   & 0   & 1    & 1   \\
    $d$   & 0   & 0   & 0   & 0    & 0   \\ \hline
    $\sum _{i} S_{ij}$   & 1   & 1   & 0 & 1.5  &     \\
    $\tilde{S} _{j} $ & 1   & 1.5 & 1 & 0  &     \\
    $R_{j}$   & 3 & 2   & 2   & 1  &     \\ \hline
    \end{tabularx}
\end{table}

\begin{table}
\centering
\caption{Flow $q_{ij}$}
\label{tab:nodemod2}
\centering
    \begin{tabularx}{\textwidth}{Y | YYYY | Y}
    \hline
    $q_{ij}$   & $a'$  & $b'$  & $c'$  & $d'$ & $q_{i}$   \\ \hline
    $a$   & 0   & 0.5   & 0   & 0     & 0.5   \\
    $b$   & 1   & 0   & 0   & 0.5  & 1.5 \\
    $c$   & 0   & 0   & 0   & 0.5  & 0.5   \\
    $d$   & 0   & 0   & 0   & 0  & 0   \\ \hline
    $\sum _{i} q_{ij}$   & 1   & 0.5   & 0 & 1  &     \\
    $\sum _{i} q_{ij} + \tilde{S} _{j}$   & 2 & 2   & 1   & 1  &     \\ \hline
    \end{tabularx}
\end{table}

In Table \ref{tab:nodemod1}, the first column represents the demand toward link $a'$ which is originated from only link $b$ as $S_{ba'}$. From Equation \eqref{subeq:nodemodel_c}, the demand constraint limits the flow $q_{ba'} \leq S_{b}$ which must be less than 1.5. From the supply constraint Equation \eqref{subeq:nodemodel_d}, we can rewrite the equation to determine the total flow as $q_{ba'} \leq R_{a'} - \tilde{S}_{a'}$. Since the supply constraint $R_{a'}$ is equal to 3 and the {\color{black}demand on the same bidirectional} link $a'$ can be interpreted as $\tilde{S}_{a'} = S_{a} = 1$, the total flow $q_{a'}$ will be no greater than 2. From the Equation \eqref{subeq:nodemodel_e}, we can rewrite the partial flow as  $q_{ba'}  = \frac{S_{ij} q_{i} }{S_{i}} = \frac{S_{ba'} q_{b} }{S_{b}} = \frac{(1)(1.5)}{1.5} = 1$. This shows that the partial flow $q_{ba'}$ is limited by the demand constraint.

The second column represents the demand toward link $b'$ which is originated from only link $a$ as $S_{ab'}$. From Equation \eqref{subeq:nodemodel_c}, the demand constraint limits the flow $q_{ab'} \leq S_{a}$ which must be less than 1.
From the supply constraint in Equation \eqref{subeq:nodemodel_d}, the total flow toward link $b'$ can be expressed as $q_{ab'} \leq R_{b'} - \tilde{S}_{b'}$. Since the supply constraint $R_{b'}$ is equal to 2 and the total demand on the same bidirectional flow can be expressed as $\tilde{S}_{b'} = S_{b} = S_{ba'} + S_{bd'} = 1+0.5 = 1.5$, the total flow $q_{b'}$ is no greater than 0.5. From the Equation \eqref{subeq:nodemodel_e}, the partial flow can be expressed as $q_{ab'}  = \frac{S_{ij} q_{i} }{S_{i}} = \frac{S_{ab'} q_{a} }{S_{a}} = \frac{(1)(0.5)}{1} = 0.5$. This shows that the partial flow $q_{ab'}$ is limited by the supply constraint.

The fourth column represents demand toward link $d'$ which are originated from link $b$ and $c$ as $S_{bd'}$ and $S_{cd'}$, respectively. From Equation \eqref{subeq:nodemodel_c}, the demand constraint limits the total flow $q_{d'} \leq S_{bd'} + S_{cd'}$ which must be less than 1.5. From the supply constraint in Equation \eqref{subeq:nodemodel_d}, the total flow toward link $d'$ can be expressed as $q_{d'} \leq R_{d'} - \tilde{S}_{d'}$. Since the supply constraint $R_{b'}$ is equal to 1 and the total demand on the same bidirectional flow can be expressed as $\tilde{S}_{d'} = S_{d} = 0$, the total flow $q_{d'}$ is no greater than 1. {\color{black}From the Equation \eqref{subeq:nodemodel_e}, the partial flow can be expressed as $q_{bd'}  = \frac{S_{bd'} q_{b}}{S_{b}} = \frac{(0.5)(1.5)}{1.5} = 0.5$.} Similarly, $q_{cd'} = \frac{S_{cd'} q_{c} }{S_{c}} = \frac{(1)(0.5)}{1}$ = 0.5. The partial flow $q_{cd'}$ is limited by the supply constraint. However, the partial flow $q_{bd'}$ is limited by both the demand and supply constraints.

{\color{black}
Both the route choice model and DNL model must work together successively. The route choice model provides time-varying turning fractions specifically on each destination. Then the DNL utilizes turning fractions to simulate the propagation of pedestrians in the network under changing traffic states throughout the simulation. The first iteration of the route choice model performs on the empty network, so no pedestrian has been loaded into the network. The route choice model is based on a free flow route travel time directly proportional to distance. In this case, undertaken paths are simply the shortest paths. However, from the second iteration, the network has pedestrians still in the middle of their journey on the network. In this case, the route choice model adjusts to consider pedestrian traffic delays from bidirectional traffic conditions following pVDFs.}

\section{Numerical Experiments} \label{sec:example}
{\color{black}
This section provides numerical experiments including two hypothetical networks and one real-world large-scale network to demonstrate the applicability and validity of the proposed pedestrian DTA modeling framework. Table \ref{tab:numerical_overview} provides an overview summary of the settings for each of the numerical experiment scenarios.
\begin{table}
\centering
\caption{Summary of the numerical experiment scenarios}
\label{tab:numerical_overview}

\begin{tabularx}{\textwidth}{Y|YYYYYYY}
\hline
  Scenario & 1 & 2 & 3$^*$ & 4 & 5 & 6 & Sydney \\ \hline
Nodes & 9 & 9 & 9 & 10 & 10 & 10 & 3,341 \\
Links & 24 & 24 & 24 & 18 & 18 & 18 & 19,612 \\
OD pairs & 1 & 2 & 1 & 1 & 2 & 2 & 21,762 \\\hline
\end{tabularx}
\footnotesize{$^*$ Significantly increase travel time of Link 7-8  after 20 seconds }
\end{table}

To demonstrate that the obtained flow assignments are converging to DUE, we estimate the relative gap \citep{chiu2011dynamic} as expressed below

\begin{equation} \color{black}
rel_{gap} = \frac{\sum _{k \in K} \sum _{r \in N} \sum _{s \in N} \sum _{p \in \Pi _{rs}} f_{p,k} ^{rs} u_{p,k}^{rs} - \sum _{k \in K} \sum _{r \in N} \sum _{s \in N} q_{rs}^{k} \hat{u} _{p,k} ^{rs}} {\sum _{k \in K} \sum _{r \in N} \sum _{s \in N} q_{rs}^{k} \hat{u} _{p,k} ^{rs}}
\label{eq:rgap}
\end{equation} 
where $\hat{u} _{p,k} ^{rs}$ denotes the shortest route travel time of path $p$ at departure time $k$ between origin $r$ and destination $s$.}

\subsection{Small grid network}
\begin{figure}
    \centering
    \subfloat[]{
        \includegraphics[width=0.3\linewidth]{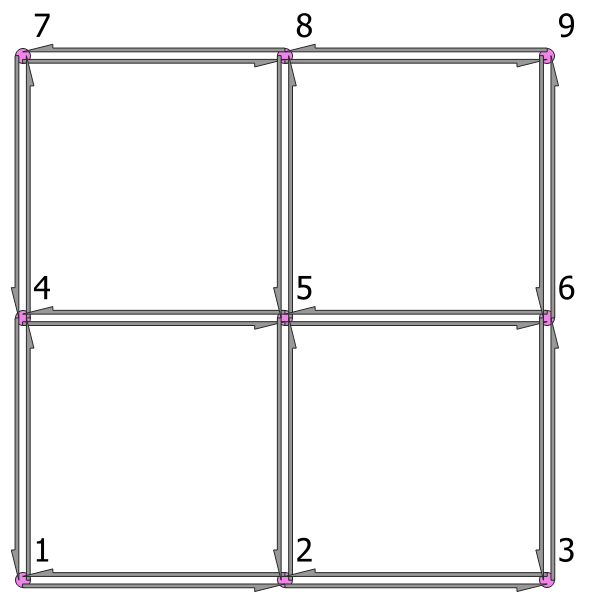}
        \label{subfig:network1_overview}} 
    \subfloat[]{
        \includegraphics[width=0.31\linewidth]{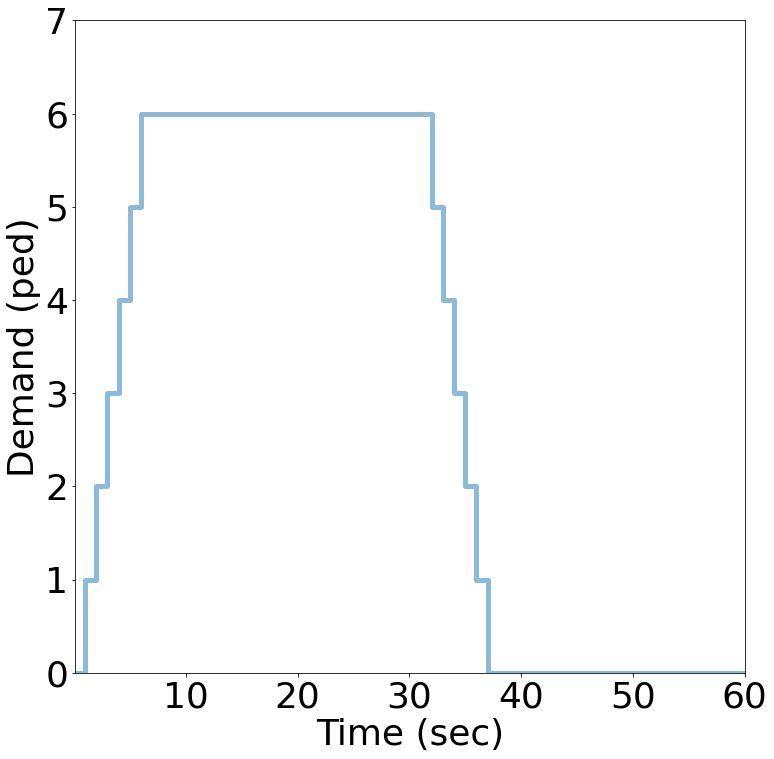}
        \label{subfig:case1_demandprofile}} 
    \quad 
    \subfloat[]{
    	\includegraphics[width=0.31\linewidth]{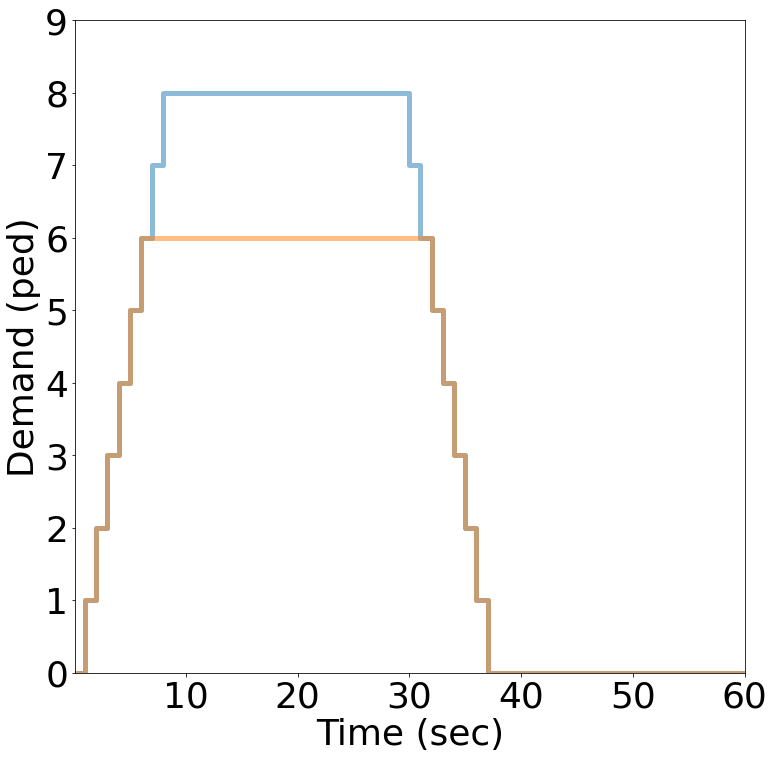}
    	\label{subfig:case2_demandprofile}}

    \caption{Small network demonstration: (a) the grid network structure. Pink circles represent nodes and gray lines with arrows represent directional links. (b) The demand profile for scenarios 1 and 3 with unidirectional flows. The black line represents demand from node 1 to node 9. (c) The demand profile for scenario 2 with bidirectional flows. The black line represents demand from node 1 to node 9 and the orange line represents demand from node 8 to node 4.  }
    \label{fig:network1_overview}
\end{figure}

\begin{figure}
\centering
\includegraphics[width=0.6\linewidth]{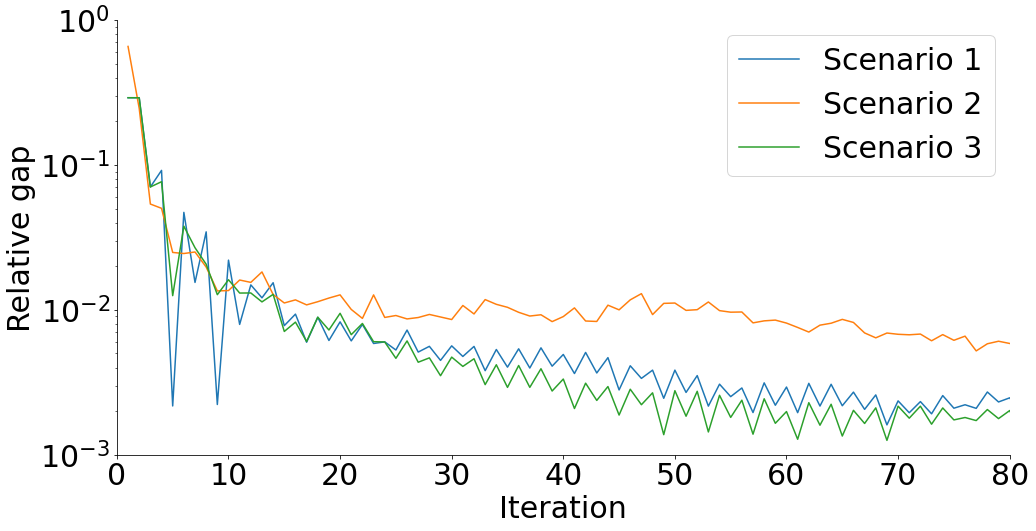}
\caption{The convergence pattern in scenario 1, 2, and 3.}
\label{fig:rgapgrid}
\end{figure}

{\color{black}We first apply the proposed DTA framework on a small walking network consisting of 9 nodes and 24 links with 12 pairs of bidirectional links for demonstration purposes. Each link in the network is 4m wide and 2m long as shown in Figure \ref{subfig:network1_overview}. We perform three hypothetical scenarios: Scenario 1 includes pedestrian travel demand for a single OD pair only, going from node 1 to node 9 as shown in Figure \ref{subfig:case1_demandprofile}. Scenario 2 includes pedestrian travel demand for two OD pairs, from node 1 to node 9 and from node 8 to node 4 as shown in Figure \ref{subfig:case2_demandprofile}. Scenario 3 includes the same travel demand as in scenario 1 as shown in Figure \ref{subfig:case1_demandprofile}; however after 20 seconds, a large travel cost is imposed on link 4-7 until the end of the simulation, creating a hypothetical link closure in the network. Figure \ref{fig:rgapgrid} reveals the convergence pattern in scenarios 1, 2 and 3. The convergence of the bidirectional traffic in scenario 2 is slower than the unidirectional traffic in scenario 1 and 3.}

\begin{figure}
    \centering
    \subfloat[]{
        \includegraphics[width=0.22\linewidth]{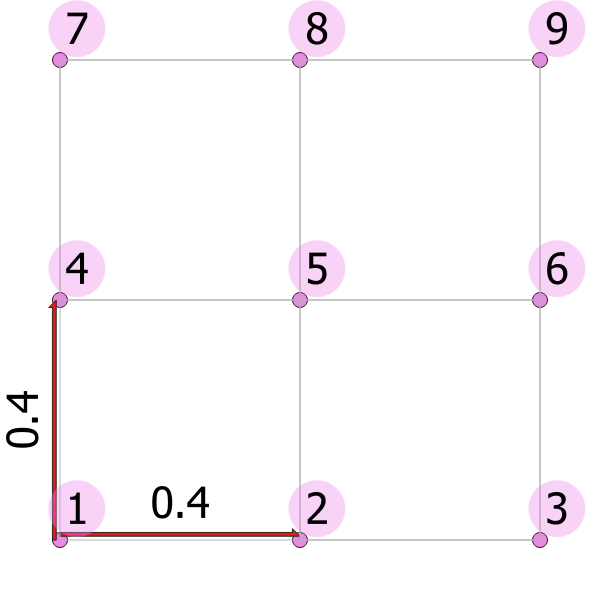}
        \label{subfig:case1_snapshot1}} 
    \hspace{0.01pt}
    \subfloat[]{
    	\includegraphics[width=0.22\linewidth]{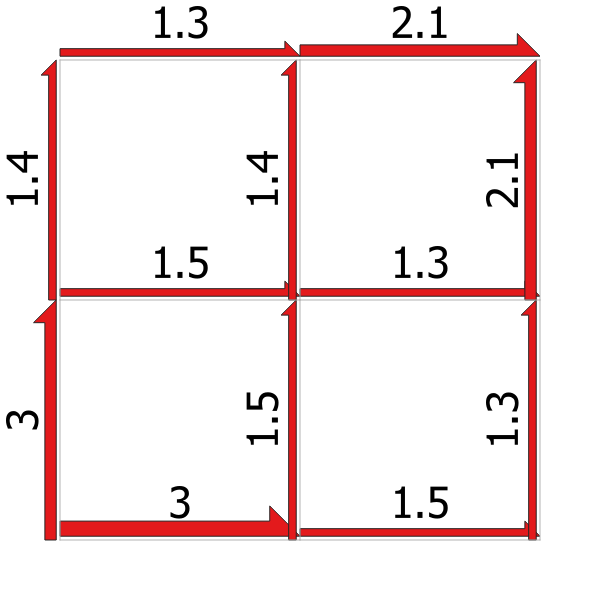}
    	\label{subfig:case1_snapshot2}} 
	\hspace{0.01pt}
	\subfloat[]{
    	\includegraphics[width=0.22\linewidth]{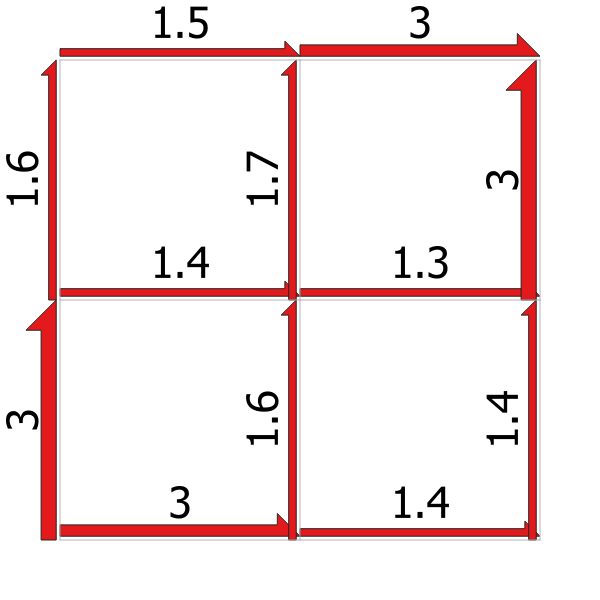}
    	\label{subfig:case1_snapshot3}} 
	\hspace{0.01pt}
	\subfloat[]{
    	\includegraphics[width=0.22\linewidth]{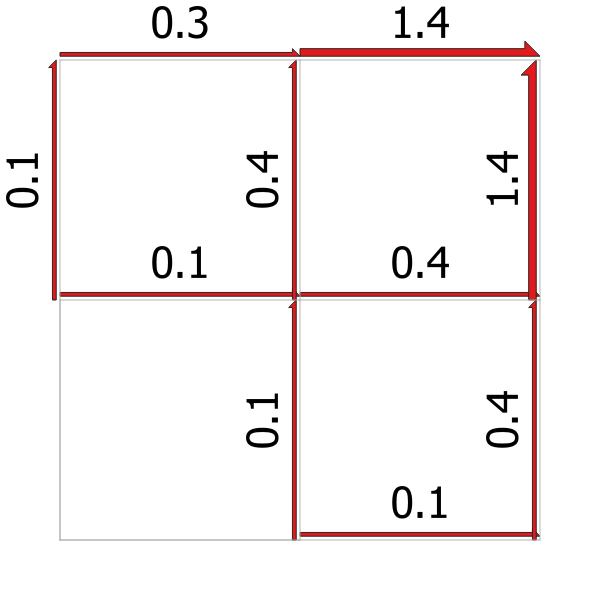}
    	\label{subfig:case1_snapshot4}} 
    
    \caption{Simulation outcomes for scenario 1 over time:  (a) t = 3 sec (b) t =10 sec (c) t = 30 sec (d) t = 40 sec. The arrow thickness and the indicated numbers represent link outflows ranging 0-4 ped/sec.}
    \label{fig:case1_snapshot}
\end{figure}

\begin{figure}
    \centering
    \subfloat[]{
        \includegraphics[width=0.22\linewidth]{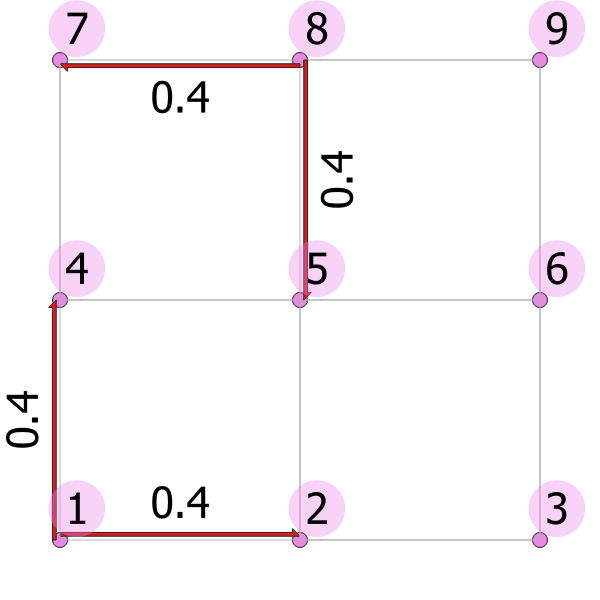}
        \label{subfig:case2_snapshot1}} 
    \hspace{0.01pt}
    \subfloat[]{
    	\includegraphics[width=0.22\linewidth]{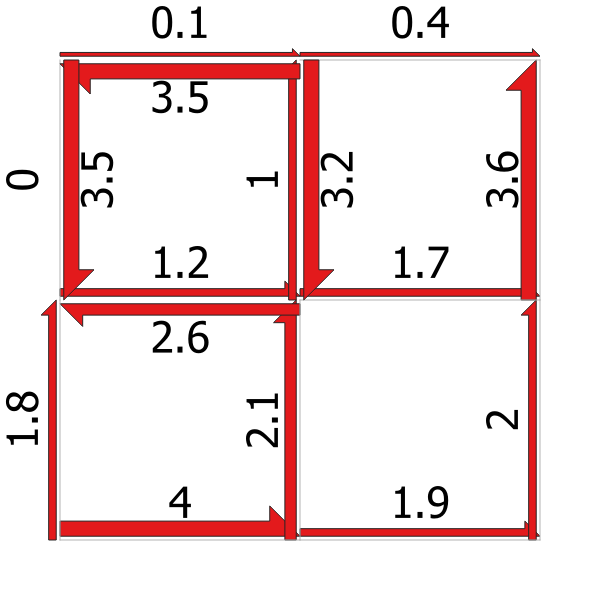}
    	\label{subfig:case2_snapshot2}} 
	\hspace{0.01pt}
	\subfloat[]{
    	\includegraphics[width=0.22\linewidth]{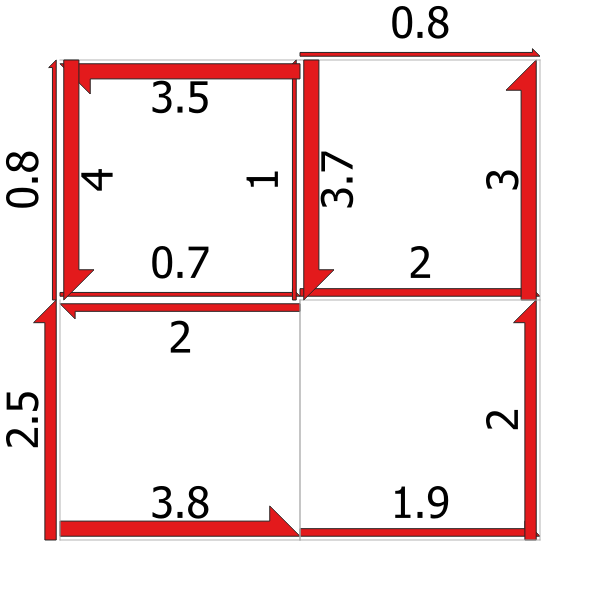}
    	\label{subfig:case2_snapshot3}} 
	\hspace{0.01pt}
	\subfloat[]{
    	\includegraphics[width=0.22\linewidth]{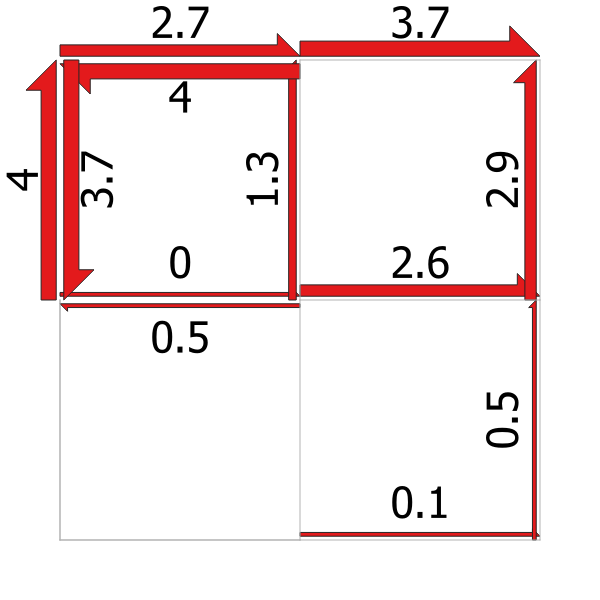}
    	\label{subfig:case2_snapshot4}} 
    
    \caption{Simulation outcomes for scenario 2 over time:  (a) t = 3 sec (b) t =10 sec (c) t = 30 sec (d) t = 40 sec. The arrow thickness and the indicated numbers represent link outflows ranging 0-4 ped/sec.}
    \label{fig:case2_snapshot}
\end{figure}

\begin{figure}
    \centering
    \subfloat[]{
        \includegraphics[width=0.22\linewidth]{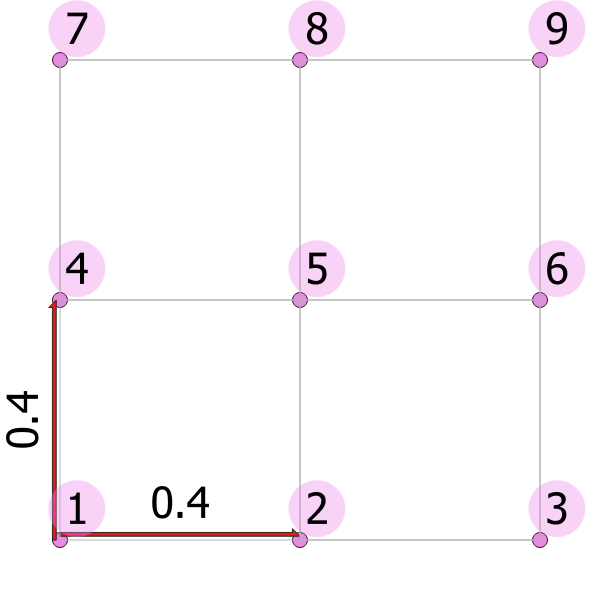}
        \label{subfig:case3_snapshot1}} 
    \hspace{0.01pt}
    \subfloat[]{
    	\includegraphics[width=0.22\linewidth]{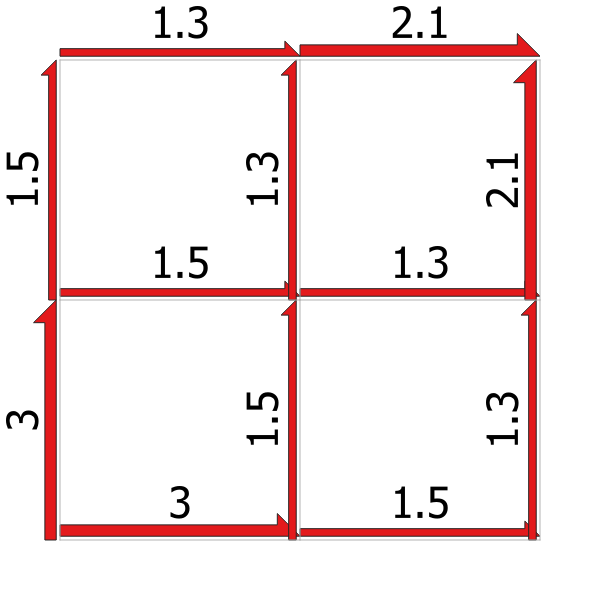}
    	\label{subfig:case3_snapshot2}} 
	\hspace{0.01pt}
	\subfloat[]{
    	\includegraphics[width=0.22\linewidth]{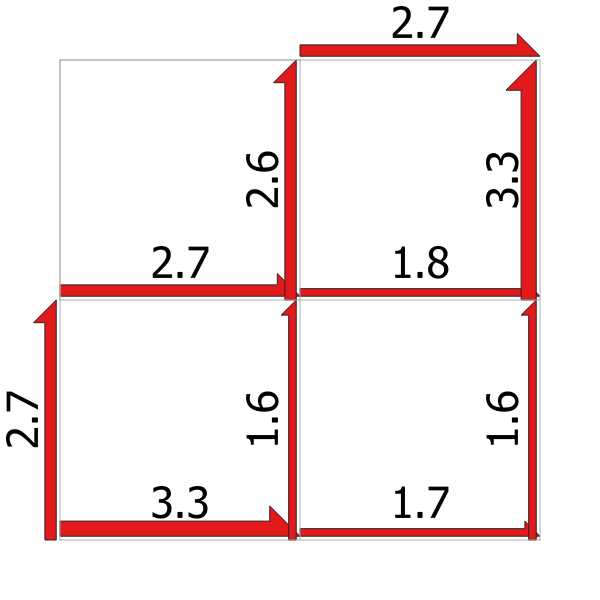}
    	\label{subfig:case3_snapshot3}} 
	\hspace{0.01pt}
	\subfloat[]{
    	\includegraphics[width=0.22\linewidth]{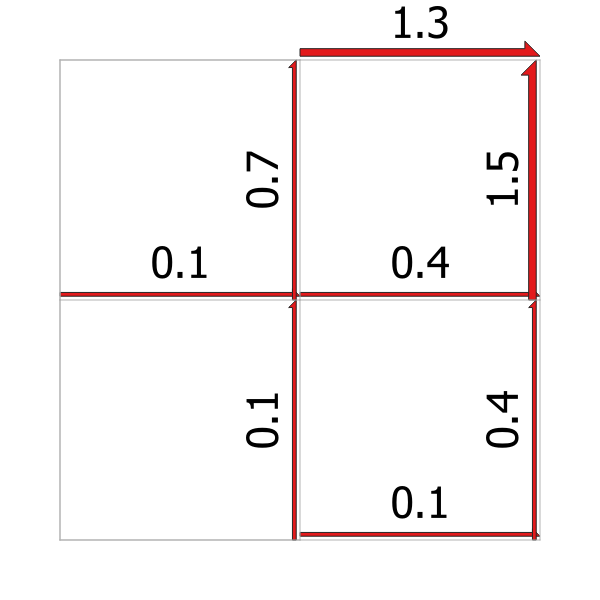}
    	\label{subfig:case3_snapshot4}} 
    
    \caption{Simulation outcomes for scenario 3 over time:  (a) t = 3 sec (b) t =10 sec (c) t = 30 sec (d) t = 40 sec. The arrow thickness and the indicated numbers represent link outflows ranging 0-4 ped/sec.}
    \label{fig:case3_snapshot}
\end{figure}

\begin{figure}
    \centering
    \subfloat[]{
        \includegraphics[width=0.25\linewidth]{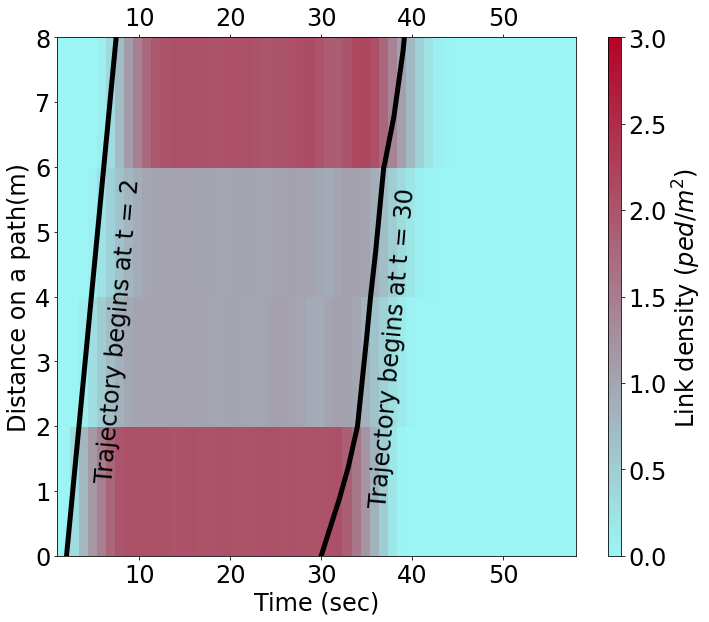}
        \label{subfig:case1_timespacedensity}} 
    \hspace{0.01pt}
    \subfloat[]{
    	\includegraphics[width=0.25\linewidth]{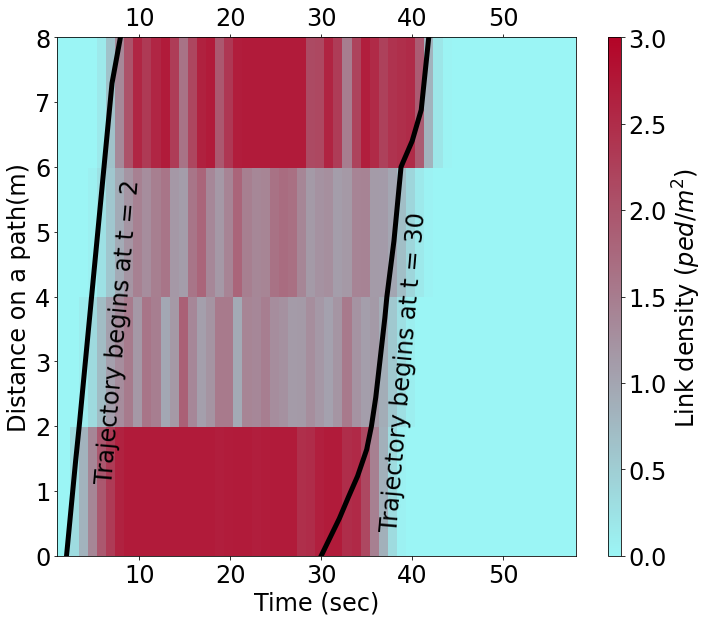}
    	\label{subfig:case2_timespacedensity}} 
	\hspace{0.01pt}
	\subfloat[]{
    	\includegraphics[width=0.25\linewidth]{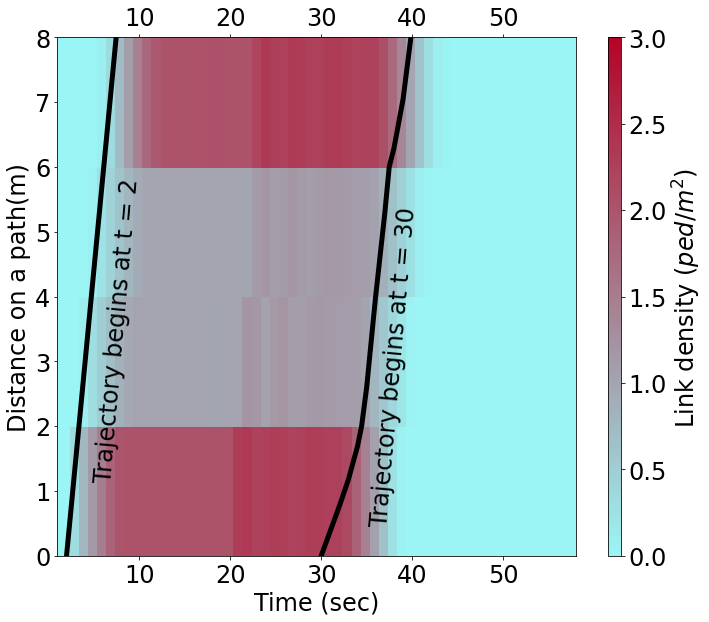}
    	\label{subfig:case3_timespacedensity}}

    \caption{Simulation outcomes for path 1-2-3-6-9 for (a) scenario 1, (b) scenario 2 and (c) scenario 3. Time-space diagrams of link density are shown in $ped/m^{2}$. Grey lines represent examplee trajectories entering the network at t = 2 sec and t = 30 seconds. }
    \label{fig:case1_timespacediagram}
\end{figure}

We treat scenario 1 as the unidirectional baseline in which high pedestrian densities are likely to occur on links emanating from the origin node 1 and going to the destination node 9. All possible paths between the single OD pair inevitably include either link 1-4 or link 1-2 at the beginning and link 6-9 and link 8-9 at the end. All other links in the network will naturally end up having lower densities since flows will be distributed uniformly across the network in unidirectional streams. {\color{black} In scenario 1, flows are distributed evenly in the network as shown in Figure \ref{fig:case1_snapshot}. Link densities on path 1-2-3-6-9 are also shown in Figure \ref{subfig:case1_timespacedensity}. In scenario 2, the network is relatively more congested at the top part (e.g., node 4,5,7, and 8) as shown in Figure \ref{fig:case2_snapshot}. Link densities on path 1-2-3-6-9 are higher than scenario 1 as shown in Figure \ref{subfig:case2_timespacedensity}. The pedestrian route choice patterns in scenario 3 are similar to scenario 1 prior to $t=20$ seconds as shown in Figure \ref{fig:case3_snapshot} and \ref{subfig:case3_timespacedensity}. Once the link closure at $t=20$ seconds is applied, link densities on path 1-2-3-6-9 increase. If a pedestrian begins a trip at $t=30$ seconds, the earliest trip completion time occurs in scenario 1, followed by scenario 3, and lastly scenario 2 as shown in Figure \ref{fig:case1_timespacediagram}. The network in scenario 2 has the highest congestion level compared to the other two scenarios. }

\begin{figure}
    \centering
    \includegraphics[width=0.6\linewidth]{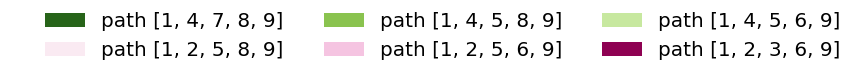}\\
    \subfloat[]{
        \includegraphics[width=0.25\linewidth]{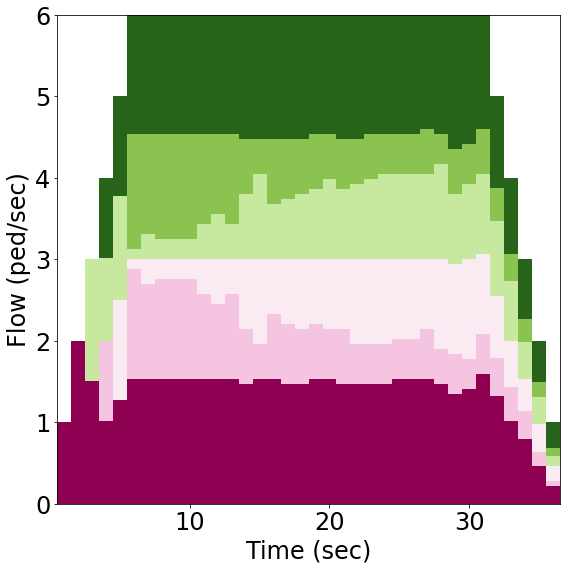}
        \label{subfig:case1_pathflow}} 
    \hspace{0.01pt}
    \subfloat[]{
    	\includegraphics[width=0.25\linewidth]{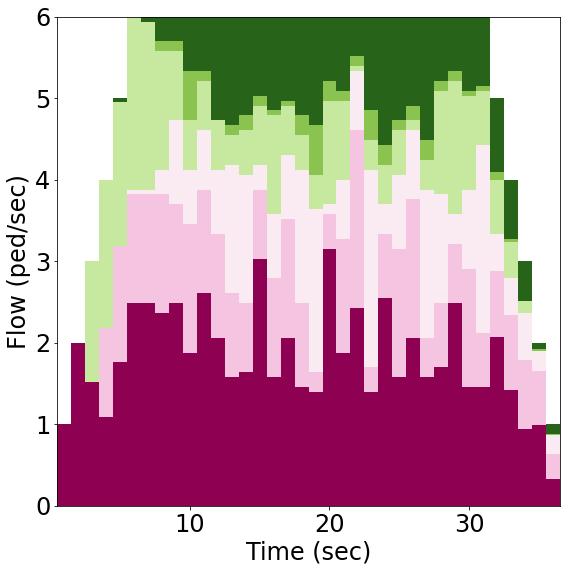}
    	\label{subfig:case2_pathflow}} 
	\hspace{0.01pt}
	\subfloat[]{
    	\includegraphics[width=0.25\linewidth]{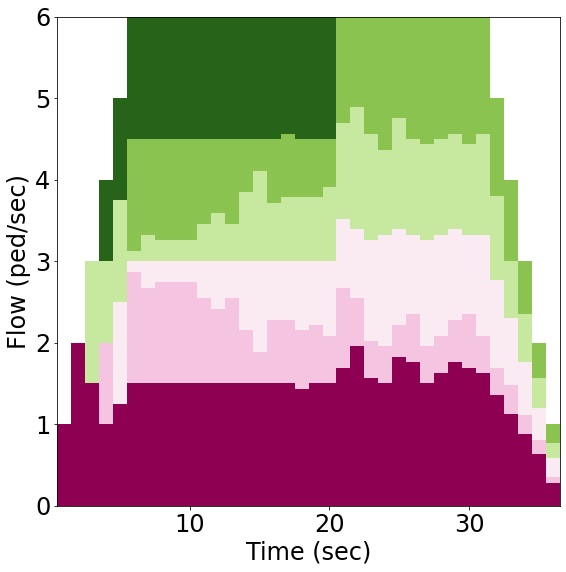}
    	\label{subfig:case3_pathflow}} 
    
    \caption{Distribution of path flows for (a) scenario 1, (b) scenario 2 and (c) scenario 3.}
    \label{fig:network_1_pathflow}
\end{figure}

Overall, there are six possible paths connecting node 1 and node 9. As shown in Figure \ref{fig:network_1_pathflow}, between $t$ = 8 sec and $t$ = 42 sec, all path flows add up to 6 pedestrians/second. Paths colored in green include link 1-4 while paths colored in pink include link 1-2. For scenario 1 as shown in Figure \ref{subfig:case1_pathflow}, flows on pink and green paths are roughly equal as expected since the grid network is totally symmetric and flows are unidirectional. For scenario 2, Figure \ref{subfig:case2_pathflow} shows that the total flows on pink paths are dominant compared to the green paths suggesting that pedestrians prefer to go through link 1-2 rather than link 1-4. {\color{black}This reflects the impact of congestion from bidirectional pedestrian flows on link 4-7 and link 4-5 because of the presence of pedestrian flow from the opposite direction on link 7-4 and link 5-4.} For scenario 3, Figure \ref{subfig:case3_pathflow} shows that after $t$ = 20 sec, no pedestrian takes path 1-4-7-8-9. This demonstrates the significant impact of introducing the link closure in the network and how pedestrians choose alternative routes when a part of the network is blocked.

\begin{figure}
    \centering
    \subfloat[]{
        \includegraphics[height=0.3\linewidth]{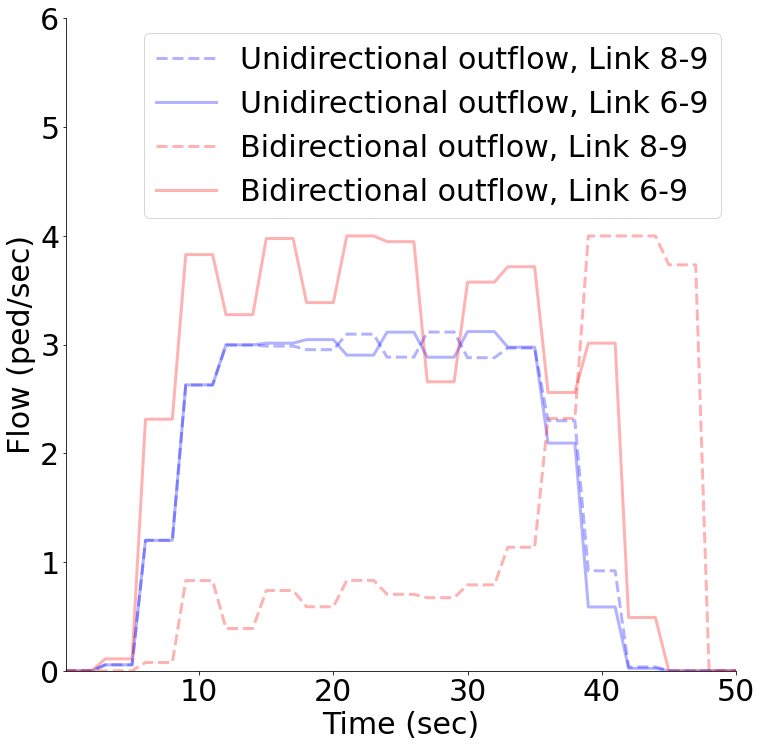}
        \label{subfig:case2_linkflow}} 
    \quad 
    \subfloat[]{
    	\includegraphics[height=0.3\linewidth]{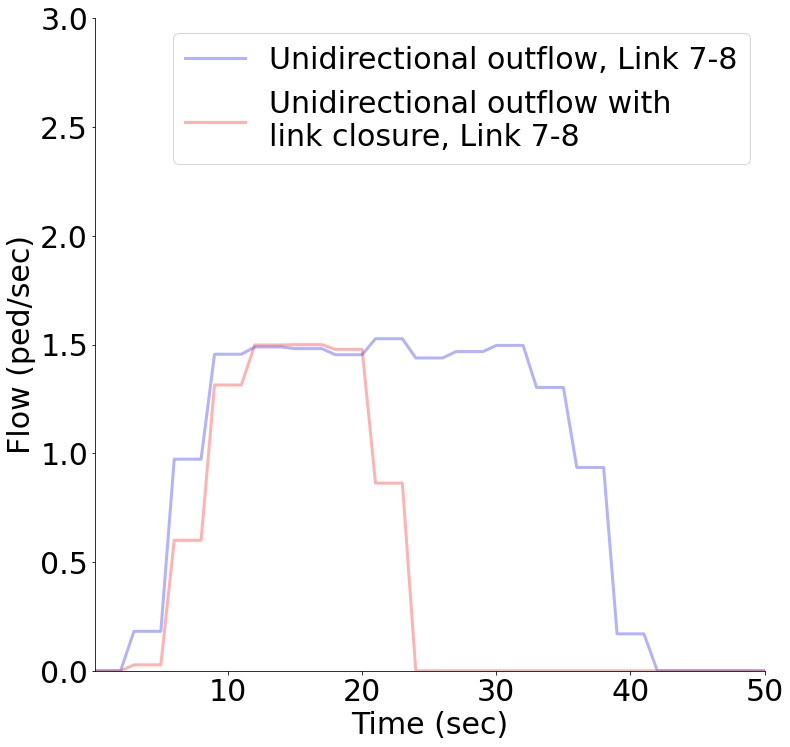}
    	\label{subfig:case3_linkflow}} 
    \caption{A comparative analysis of time-dependent link flows. Comparing scenario 1 with unidirectional flows shown in black against (a) scenario 2 with bidirectional flows shown in red for links 6-9 and 8-9 and (b) scenario 3 for link 7-8.}
    \label{fig:linkflow}
\end{figure}

As shown in Figure \ref{subfig:case2_linkflow}, in scenario 1 with unidirectional flows only, both links 8-9 and 6-9 are utilized evenly. However, in scenario 2 when bidirectional flow is present, link 6-9 almost always experiences higher flow compared with link 8-9 suggesting that pedestrians tend to avoid the congested paths due to the presence of flow from the opposite direction. Figure \ref{subfig:case3_linkflow} also shows that in scenario 3 after $t$ = 20 sec, no pedestrian uses link 7-8 due to the closure of link 4-7. While in scenario 1 without implementation of the link closure, link 7-8 was normally used by pedestrians. 
\FloatBarrier

\subsection{Long corridor}
In this section, we apply the developed DPTA framework to a long bidirectional corridor to demonstrate the impact of directionality on pedestrian flows. The corridor consists of straight walking links that are 2m long and 4m wide.

\begin{figure}
    \centering
    \subfloat[]{
    	\includegraphics[width=1\linewidth]{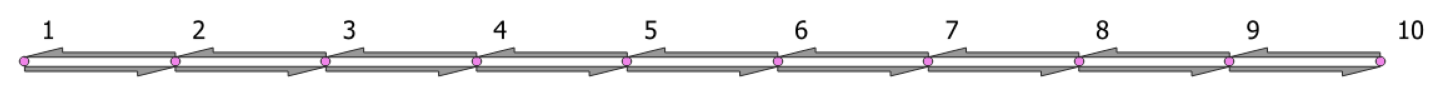}
    	\label{subfig:network2_overview}} 
	\quad
    \subfloat[]{
        \includegraphics[width=0.31\linewidth]{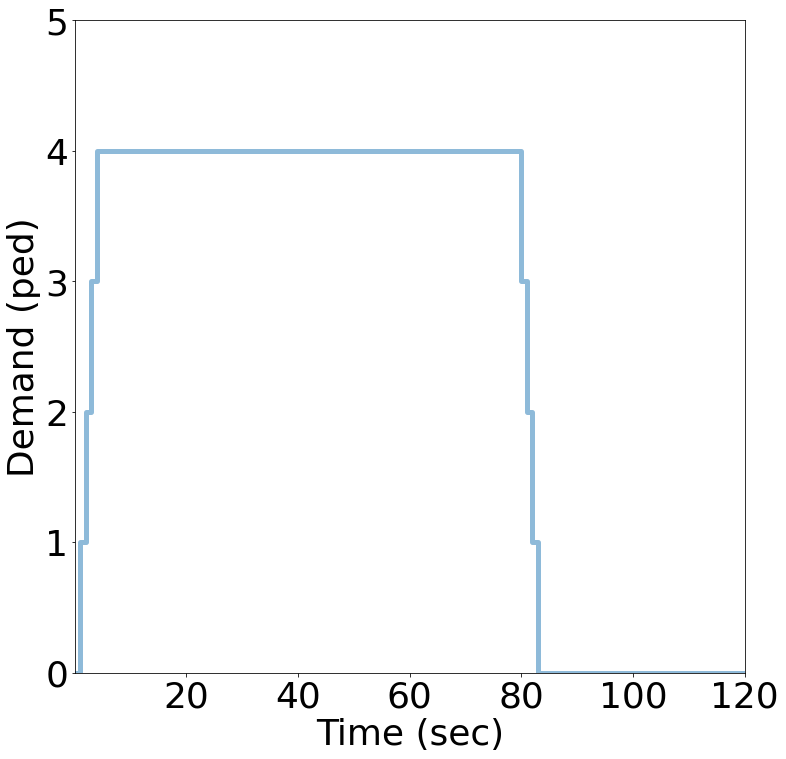}
        \label{subfig:case4_demandprofile}} 
    \hspace{0.01pt}
    \subfloat[]{
    	\includegraphics[width=0.31\linewidth]{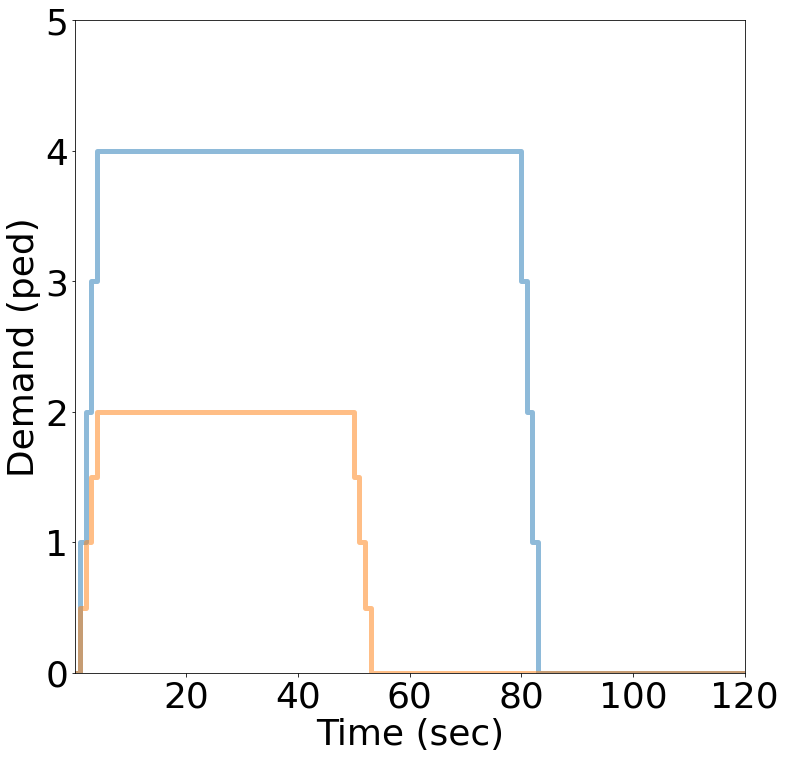}
    	\label{subfig:case5_demandprofile}} 
	\hspace{0.01pt}
	\subfloat[]{
    	\includegraphics[width=0.31\linewidth]{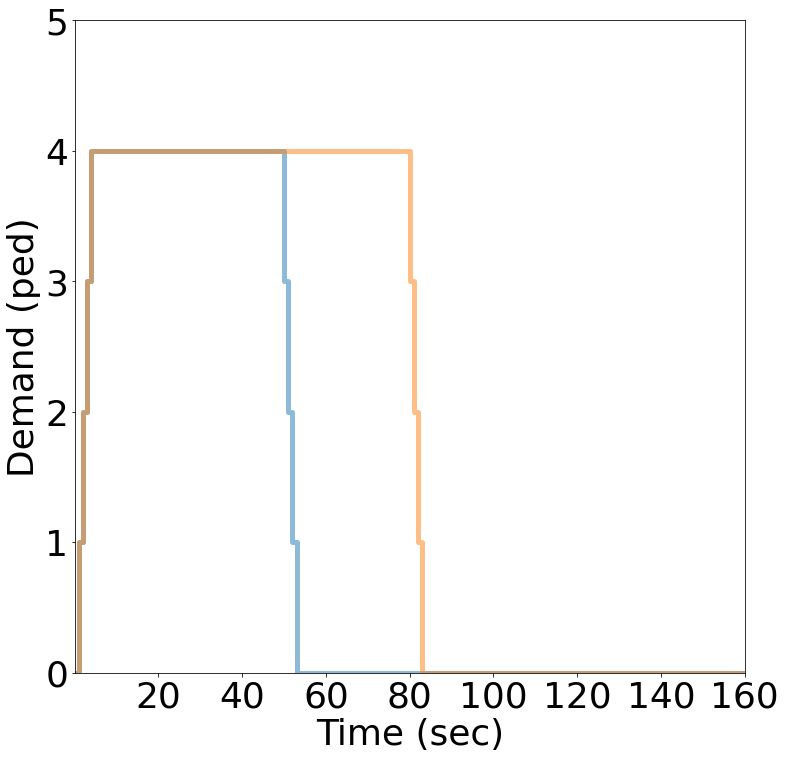}
    	\label{subfig:case6_demandprofile}} 
    
    \caption{(a) Network structure of the modeled corridor. Pink circles represent nodes and grey lines with arrows represent directional links. Simulated demand profiles for the long corridor: (b) Scenario 4 with a bottleneck, (c) Scenario 5 with unbalanced bidirectional flows, and (d) Scenario 6 with balanced bidirectional flows.{\color{black} The black line represents demand from node 1 to node 10, referred to as major stream, and the orange line represents demand from node 10 to node 1, referred as the minor stream.}}
    \label{fig:network2_demandfprofile}
\end{figure}

We perform three hypothetical scenarios numbered continuously from the previous section: In scenario 4, we decrease the width of link 9-10 to create a bottleneck and to trigger formation of a shockwave. Pedestrian demand in scenario 4 gradually increases from $t$ = 1 sec to $t$ = 4 sec reaching 4 pedestrians/second. From $t$ = 5 sec to $t$ = 80 sec, demand remains constant at 4 pedestrians/second. From $t$ = 80 sec to $t$ = 83 sec, demand gradually drops to zero and remains zero until the end of the simulation as shown in Figure \ref{subfig:case4_demandprofile}.

Scenario 5 represents a bidirectional stream with unbalanced flows consisting of pedestrian demand going from node 1 to node 10 as the major stream and demand from node 10 to node 1 in the opposite direction as the minor stream. The demand for the major direction is identical to scenario 4. However, the demand for the minor direction gradually increases between $t$ = 1 sec and $t$ = 4 sec reaching 2 pedestrians/second. From $t$ = 5 sec to $t$ = 50 sec, demand remains constant at 2 pedestrian/second. From $t$ = 50 sec to $t$ = 53 sec demand gradually drops to zero and remains zero until the end of the simulation as shown in Figure \ref{subfig:case5_demandprofile}.

Scenario 6 represents a bidirectional stream with balanced flows at the beginning gradually turning into a unidirectional stream in which demand in both directions peak at 4 pedestrians/second. The demand for the minor direction begins to drop at $t$ = 50 sec while the demand for the major direction drops later at $t$ = 80 sec as shown in Figure \ref{subfig:case6_demandprofile}.

\begin{figure}
    \centering
    \subfloat[]{
        \includegraphics[width=0.4\linewidth]{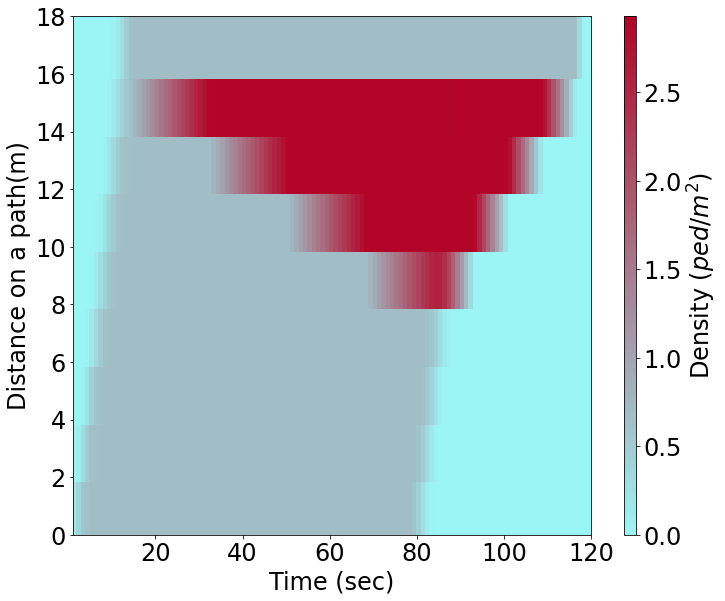}
        \label{subfig:case4_timespace_density_major}} 
    \hspace{0.01pt}
    \subfloat[]{
        \includegraphics[width=0.4\linewidth]{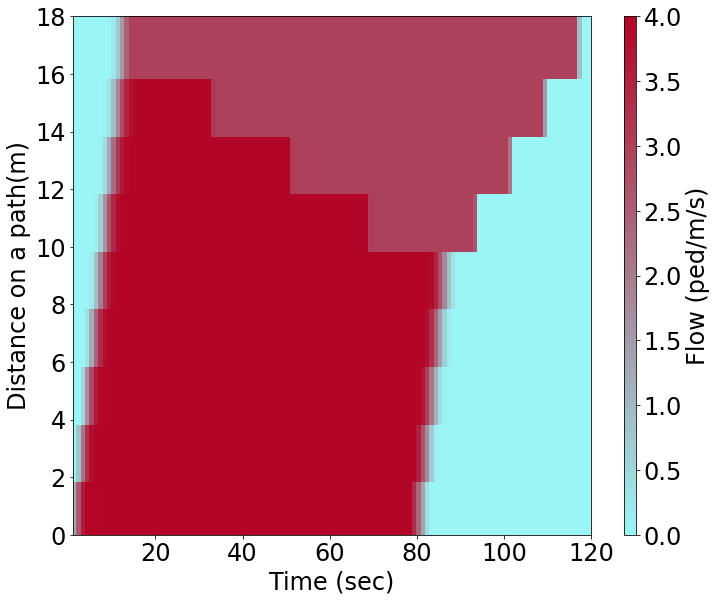}
        \label{subfig:case4_timespace_flow_major}} 
    
    \caption{Time-space diagram of scenario 4 where pedestrians travel from node 1 toward node 10: (a) density in $ped/m^{2}$ (b) flow in $ped/m/s$. }
    \label{fig:case4_timespace}
\end{figure}

{\color{black}Since the corridor has only one path and one single OD pair, the model utilizes only one iteration of route choice.} In scenario 4, pedestrians walk through the corridor with the maximum flow of 4 pedestrians/second {\color{black}from node 1 toward node 9.} Once pedestrians reach the bottleneck on link 9-10, flow begins to drop as shown in Figure \ref{subfig:case4_timespace_flow_major}. A shockwave forms on link 8-9 at distance 14-16 m and propagates backward with the speed of 0.18 m/s as shown in Figure \ref{subfig:case4_timespace_density_major}. However, after time $t$ = 85 second, as no demand loads into the network, a forward recovery shockwave forms on link 5-6 until the entire crowd exits the corridor.

\begin{figure}
    \centering
    \subfloat[]{
        \includegraphics[width=0.4\linewidth]{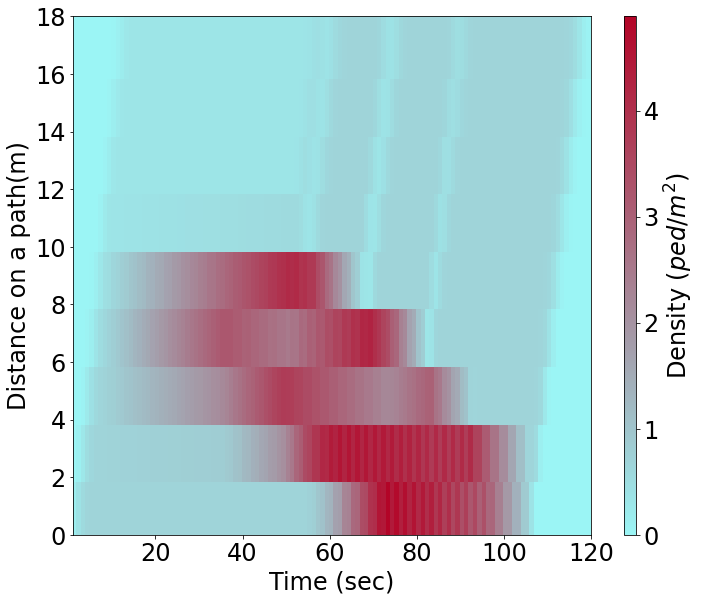}
        \label{subfig:case5_timespace_density_major}} 
    \hspace{0.01pt}
    \subfloat[]{
    	\includegraphics[width=0.4\linewidth]{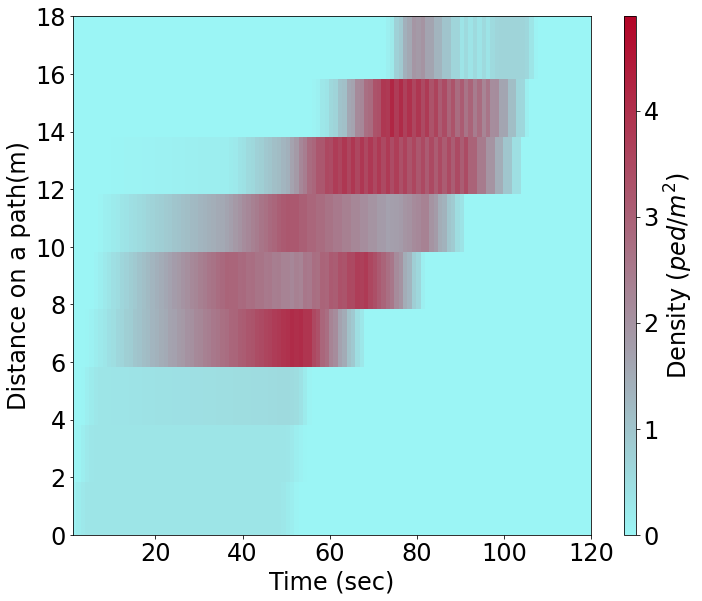}
    	\label{subfig:case5_timespace_density_minor}} 
	\hspace{0.01pt}
	   \subfloat[]{
        \includegraphics[width=0.4\linewidth]{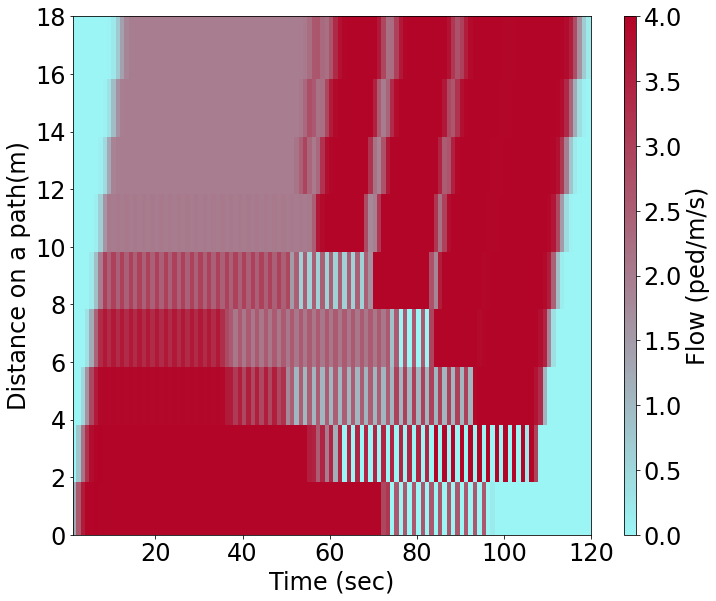}
        \label{subfig:case5_timespace_flow_major}} 
    \hspace{0.01pt}
    \subfloat[]{
    	\includegraphics[width=0.4\linewidth]{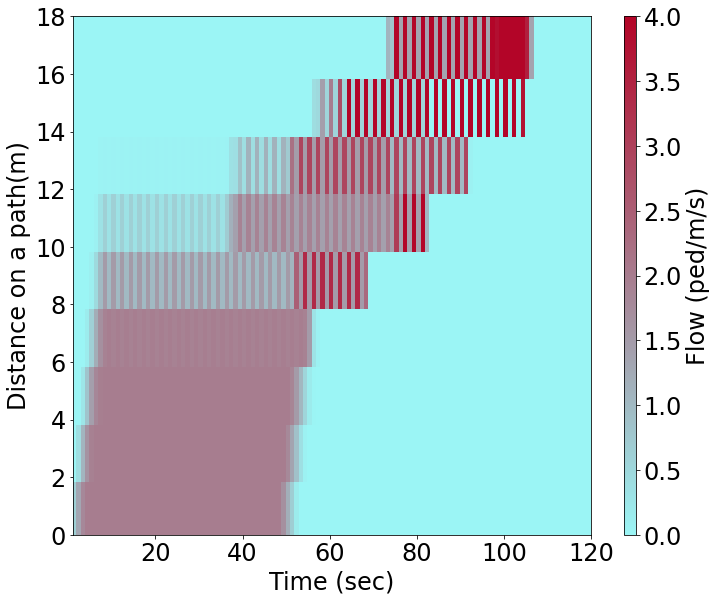}
    	\label{subfig:case5_timespace_flow_minor}} 
    \caption{Time-space diagram of scenario 5. (a) density in $ped/m^{2}$ of the major direction from node 1 toward node 10, (b) density in $ped/m^{2}$ of the minor direction from node 10 toward node 1, (c) flow in $ped/m/s$ of the major direction from node 1 toward node 10, and (d) flow in $ped/m/s$ of the minor direction from node 10 toward node 1.}
    \label{fig:case5_timespace}
\end{figure}

In scenario 5, pedestrians move through the corridor with the maximum flow of 4 pedestrians/second for the major stream and 2 pedestrians/second for the minor stream as shown in Figure \ref{subfig:case5_timespace_flow_major} and \ref{subfig:case5_timespace_flow_minor}. Once both pedestrian streams meet in the middle of the corridor at distance = 8 m, flows from both directions significantly drop. For the major direction, flow drops to 2 pedestrians/second from distance 10-18 m because the minor stream already occupies the available space to accommodate the flow of 2 pedestrians/second. The minor direction flow drops significantly at distance = 8 m as shown in Figure \ref{subfig:case5_timespace_flow_minor}. The walking speed in the minor direction also significantly reduces as pedestrians slowly make their way through as can be observed from the sharp reduction in the trajectories slope in the time-space diagram. Once the minor direction reaches node 1, it begins to block the major  stream from entering the corridor as not all the flow can load into the network as shown in Figure \ref{subfig:case5_timespace_flow_major}. When the two opposite streams meet, pedestrians take turn traversing through each node. Two shockwaves are also formed. As flow begins to drop in the major direction, a backward forming shockwave forms and propagates from distance 8 m as shown in Figure \ref{subfig:case5_timespace_density_major}. As the entire corridor is affected by the major direction, the minor direction stream exhibits a forward forming shockwave slowly penetrating through {\color{black}the major direction as shown in Figure \ref{subfig:case5_timespace_density_minor}. In microscopic pedestrian models, gridlock is a common occurrence once the bidirectional traffic condition reaches the critical density, which causes pedestrians to halt all movements \citep{muramatsu1999jamming,nagatani2009freezing,nowak2012quantitative,nowak2013cellular,tao2017cellular}. However, there is no empirical evidence showing the existence of gridlock in the literature. In reality, pedestrian movements have high degrees of freedom leading to space adaptation to avoid gridlock occurrence. To overcome this technical limitation, \citet{nowak2013cellular} and \citet{tao2017cellular} propose mechanisms to fix a gridlock. In our proposed node model, we propose a mechanism to avoid gridlock. Equation \eqref{subeq:nodemodel_d} introduces a look ahead behavior that pedestrians who want to go through a node would check in the opposite direction to ensure that flows do not exceed the receiving flow. This mechanism causes stuttering behaviors in which pedestrians from each direction take a turn going through a node to avoid gridlock. }

\begin{figure}
    \centering
    \subfloat[]{
        \includegraphics[width=0.4\linewidth]{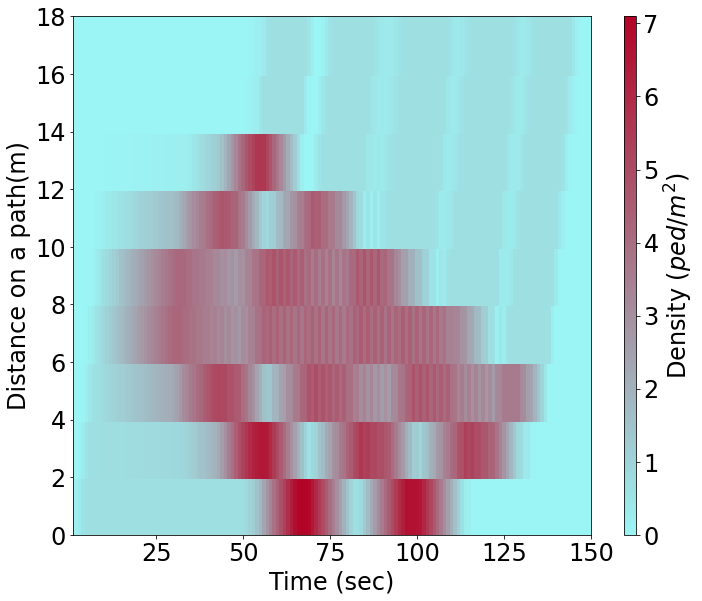}
        \label{subfig:case6_timespace_density_major}} 
    \hspace{0.01pt}
    \subfloat[]{
    	\includegraphics[width=0.4\linewidth]{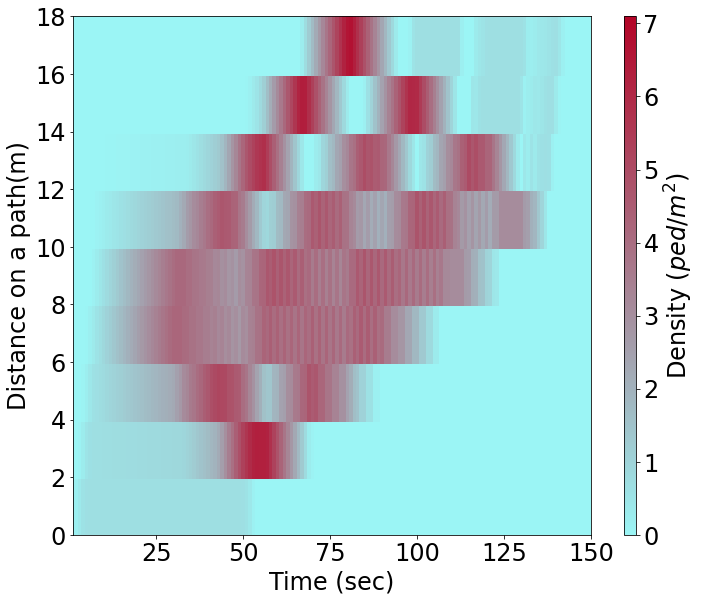}
    	\label{subfig:case6_timespace_density_minor}} 
	\hspace{0.01pt}
	   \subfloat[]{
        \includegraphics[width=0.4\linewidth]{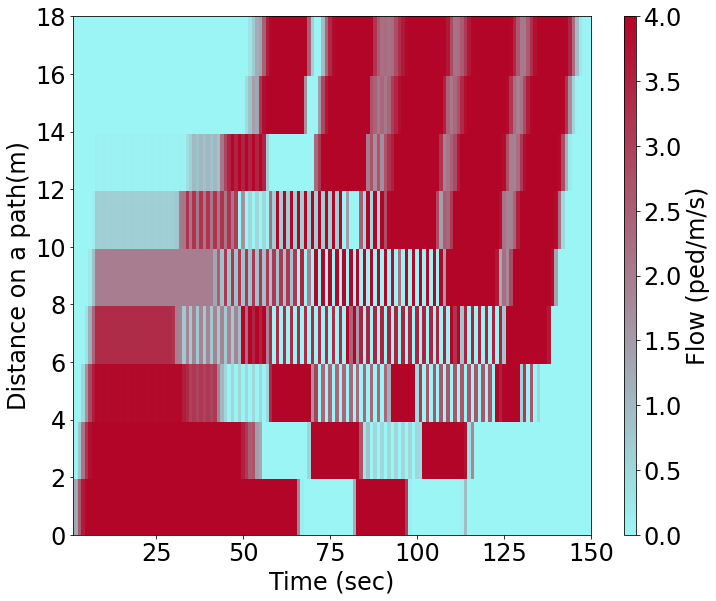}
        \label{subfig:case6_timespace_flow_major}} 
    \hspace{0.01pt}
    \subfloat[]{
    	\includegraphics[width=0.4\linewidth]{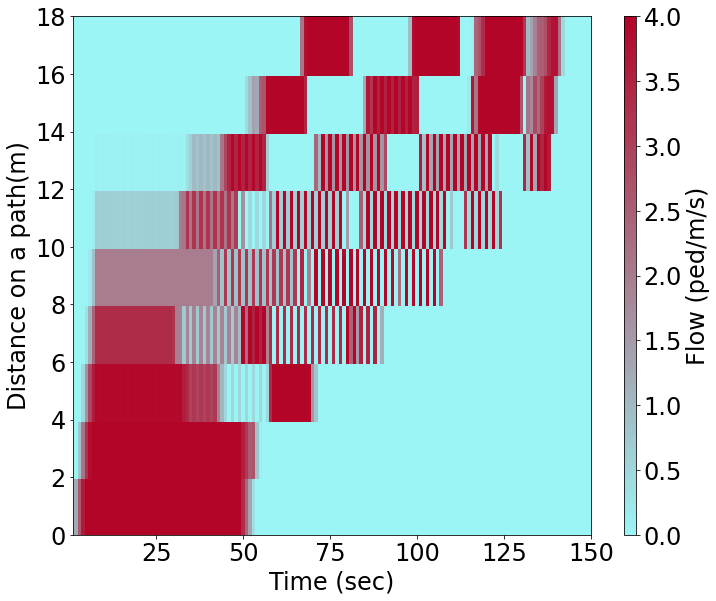}
    	\label{subfig:case6_timespace_flow_minor}} 
    \caption{\color{black}Time-space diagram of scenario 6. (a) density in $ped/m^{2}$ of the major direction from node 1 toward node 10, (b) density in $ped/m^{2}$ of the minor direction from node 10 toward node 1, (c) flow in $ped/m/s$ of the major direction from node 1 toward node 10, and (d) flow in $ped/m/s$ of the minor direction from node 10 toward node 1.}
    \label{fig:case6_timespace}
\end{figure}

In scenario 6, four shockwaves are formed in both directions as shown in Figure \ref{subfig:case6_timespace_density_major} and \ref{subfig:case6_timespace_density_minor}.{\color{black} A backward forming shockwave is formed at $t = 30$ sec that propagates upstream due to significant increase in flow and density in the major direction. A forward forming shockwave is also formed at $t = 30$ sec that propagates downstream as pedestrians from the opposite direction mix into the main direction flow. Two recovering shockwaves also occur in both directions after congestion begin to dissipate from $t = 60$ sec.} Flows in both directions significantly reduce once the two opposite streams meet at distance 8 m, then both directions take turn passing through nodes as shown in Figure \ref{subfig:case6_timespace_flow_major} and  \ref{subfig:case6_timespace_flow_minor}.
{\color{black}The simulated bidirectional long corridor successfully exhibits the formation and propagation of shockwaves when link capacities change due to the bidirectional effects or inclusion of a bottleneck. Results suggest that the proposed DPTA modeling framework reproduces reasonably realistic pedestrian dynamics compared to a static traffic assignment model.}

\FloatBarrier

\subsection{Large-scale real network application: Sydney footpath network}
{\color{black}Here, we apply the proposed DPTA framework to model the sidewalk network in Sydney CBD as a large-scale example consisting of 2,327 nodes and 8,062 links. We have estimated a total travel demand of 28,414 walking trips during a 20-minute time window in the morning peak from 21,762 OD pairs based on information obtained from local transport authorities reports \citep{transport2013} and commercially available data purchased from DSpark based on passively collected mobile phone traces. DSpark provides number of trips distinguished by travel mode between origin-destination pairs at the statistical area level 1 (SA1) aggregation. In this study, we only consider walking trips from/to the Sydney CBD area during an average weekday in October 2019 for a 20-minute interval in the morning peak. The origin-destination matrix is discretized from the SA1 level into smaller block level consistent with designed centroids in the network.} Since the focus of this study is on formulation and development of a DPTA model using LTM, calibration of the OD walking demand remains a direction for future research. {\color{black} For this numerical experiment, we only conducted 3 iterations due to the size of the network and computational limitations.}

\begin{figure}
    \centering
    \includegraphics[width=0.6\linewidth]{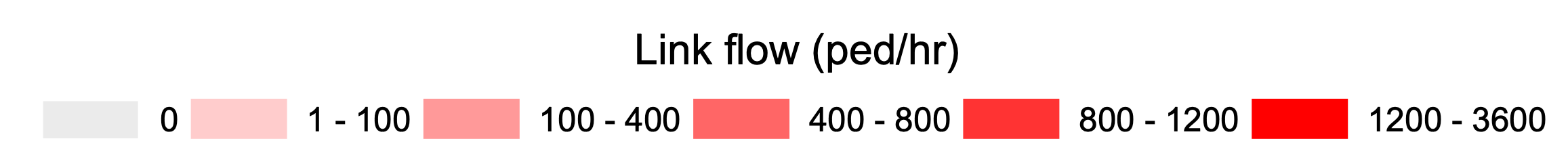}\\
    \subfloat[]{
        \includegraphics[width=0.23\linewidth]{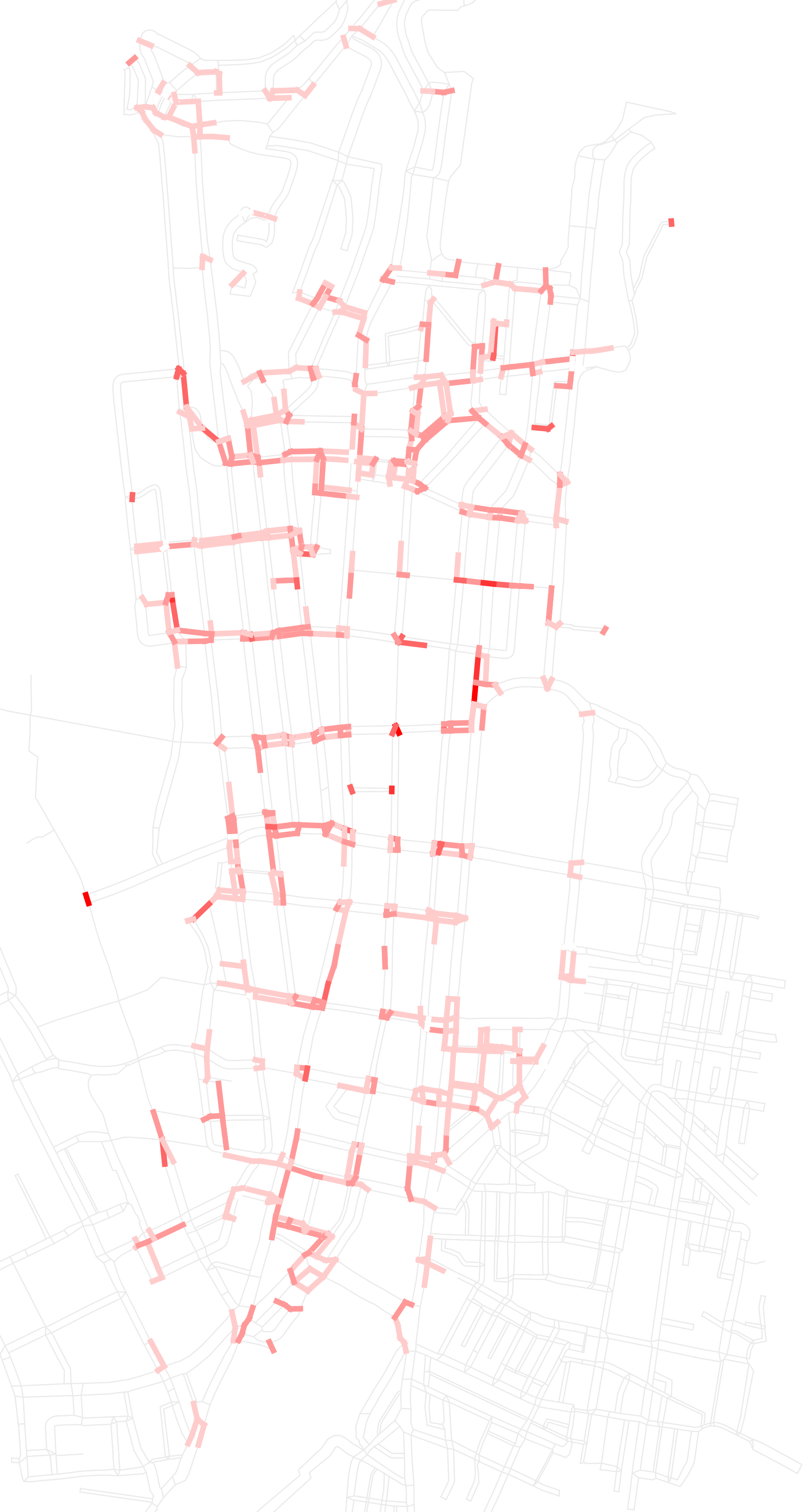}
        \label{subfig:cbd_snaphot1}} 
    \hspace{0.01pt}
    \subfloat[]{
    	\includegraphics[width=0.23\linewidth]{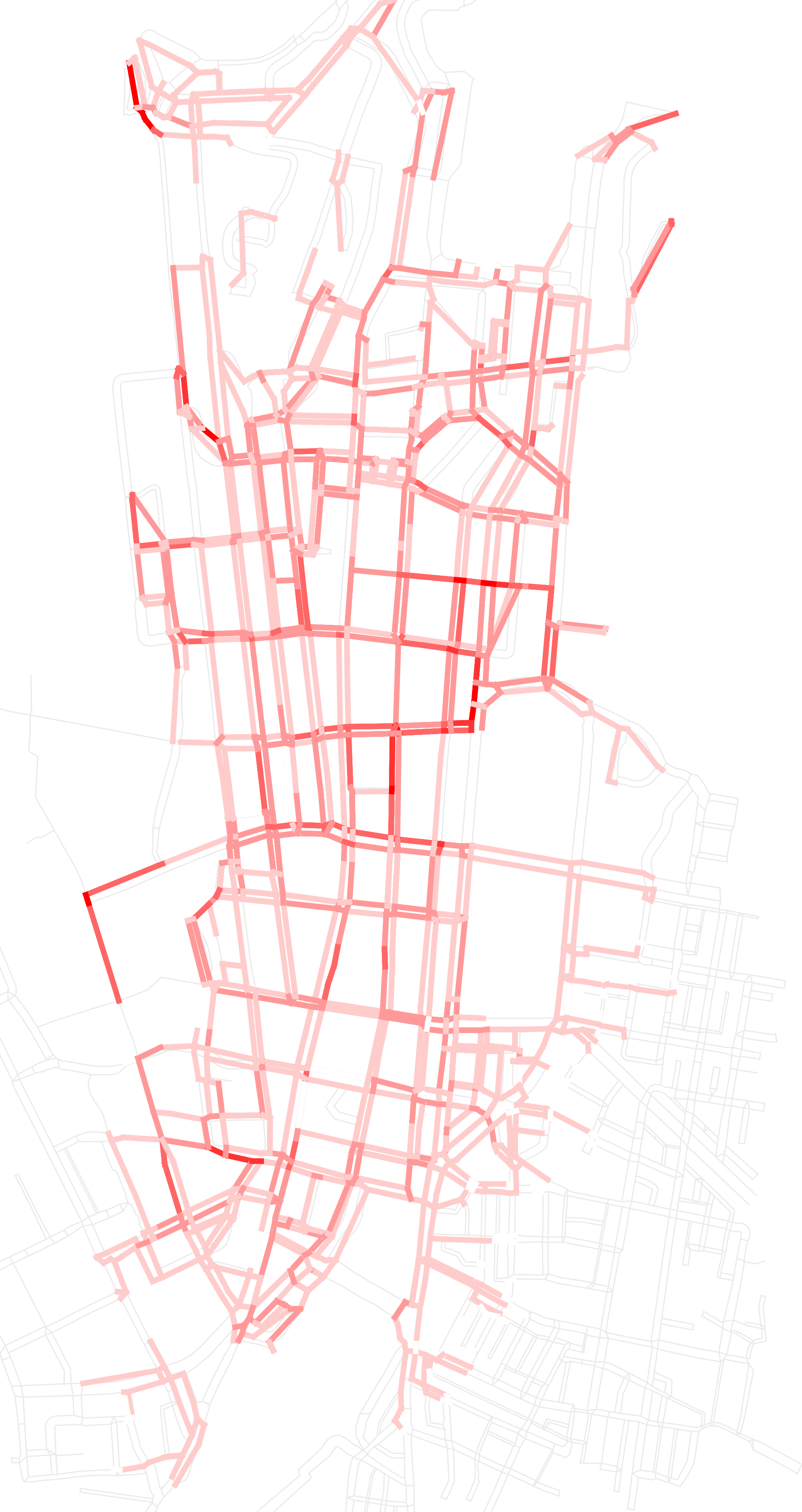}
    	\label{subfig:cbd_snapshot2}} 
	\hspace{0.01pt}
	\subfloat[]{
    	\includegraphics[width=0.23\linewidth]{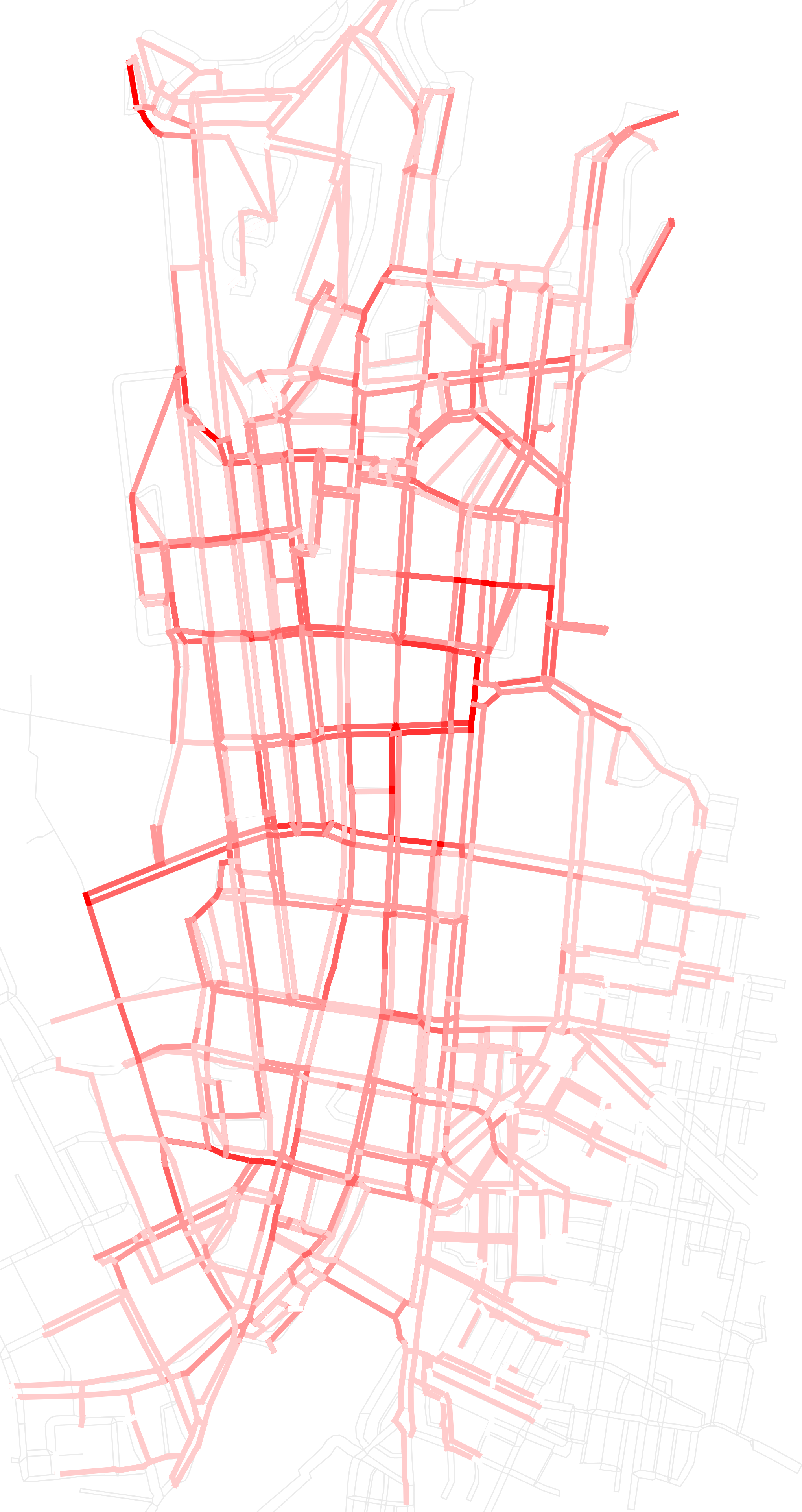}
    	\label{subfig:cbd_snapshot3}} 
	\hspace{0.01pt}
	\subfloat[]{
    	\includegraphics[width=0.23\linewidth]{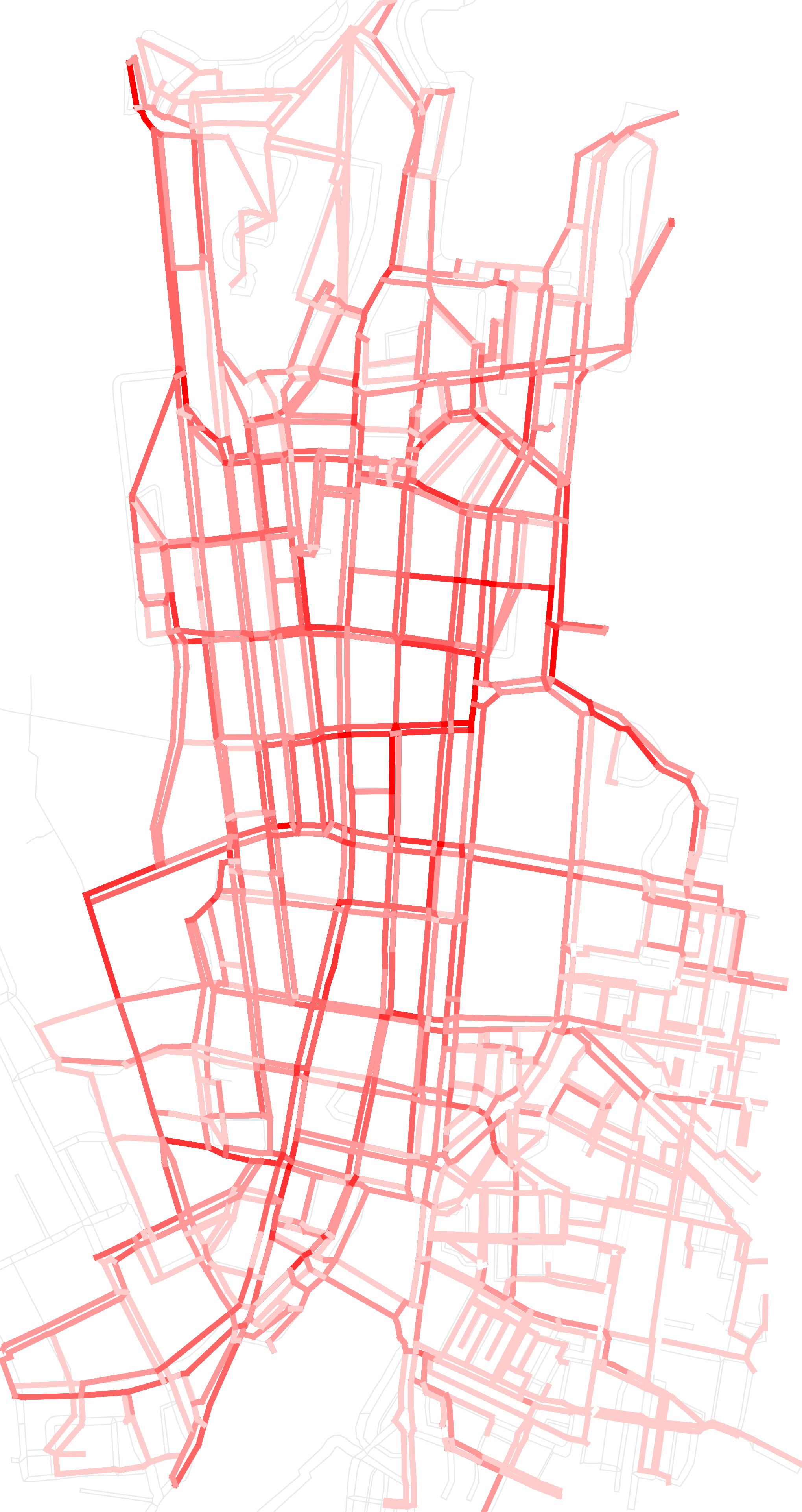}
    	\label{subfig:cbd_snapshot4}} 
    
    \caption{Simulation results of the Sydney footpath network at different times. Line thickness represents link outflows ranging from 0 - 3600 ped/hour. (a) t = 1 min (b) t = 3 min (c) t = 5 min (d) t = 14 min }
    \label{fig:cbd_snapshot}
\end{figure}

Figure \ref{fig:cbd_snapshot} shows the estimated link flows across the Sydney CBD network over time. At the beginning of the simulation, only a few links closer to public transportation stations, as the more significant origin nodes, are used as shown in Figure \ref{subfig:cbd_snaphot1} and \ref{subfig:cbd_snapshot2}. After 5 minutes, the footpath network becomes more congested as more links are being used as shown in Figure \ref{subfig:cbd_snapshot3} and \ref{subfig:cbd_snapshot4}.

{\color{black} The expected link travel time here is calculated from the pVDF as shown in Equation \eqref{eq:detvdf_asym}, while the experienced link travel time is calculated using the effective free flow speed from the FD as expressed previously in Equation \eqref{eq:freeflowspeed}. The expected route travel time in the route choice model is calculated instantaneously, while the experienced route travel time is calculated sequentially based on the entering time of each link. At the departure time, the expected or instantaneous travel time is considered as pre-trip information. However, post-trip information is necessary to determine the experienced travel time \citep{chiu2011dynamic}. We formulated an instantaneous route travel time as shown in Equation \eqref{subeq:routett_a} and an experienced route travel time as shown in Equation \eqref{subeq:routett_b} similar to formulations proposed by \citet{yildirimoglu2013experienced}. 

\begin{subequations}\label{eq:routett} 
\begin{align}
 & c_{p,I} ^{rs}(k) = \sum _{a \in A} \sum _{t \in T} c_{a,t} \delta _{a,p,k,t} ^{rs} = \sum _{i=a} ^{I} c_{i} ^{pVDF} (k)  \hspace{30pt} && \forall p \in \Pi_{rs}, \forall k \in K \label{subeq:routett_a} \\
 & c_{p,I} ^{rs} (k) = \sum _{i=a} ^{I} c_{i} ^{FD} (k+ c_{p,i} ^{rs}(k)) && \forall p \in \Pi_{rs}, \forall k \in K \label{subeq:routett_b}
\end{align}
\end{subequations}
where $c_{p,I} ^{rs}(k)$ denotes an accumulated travel time on path $p$ connecting OD pair $(r,s)$ exiting link $I$ at departure time $k$, $A$ denotes a set of links, $\Pi_{rs}$ denotes the set of paths connecting OD pair $(r,s) \in W$, $K$ denotes a set of departure time, path $p$ consists of multiple links from $a,...,I \in A$ ($a$ is the most upstream link and $I$ is the most downstream link), $c_{i} ^{pVDF}$ denotes the link travel time based on Equation \eqref{eq:detvdf_asym}, and $c_{i} ^{FD}$ denotes the link travel time derived from Equation \eqref{eq:freeflowspeed}.}

\begin{figure}
    \centering
        \includegraphics[width=\linewidth]{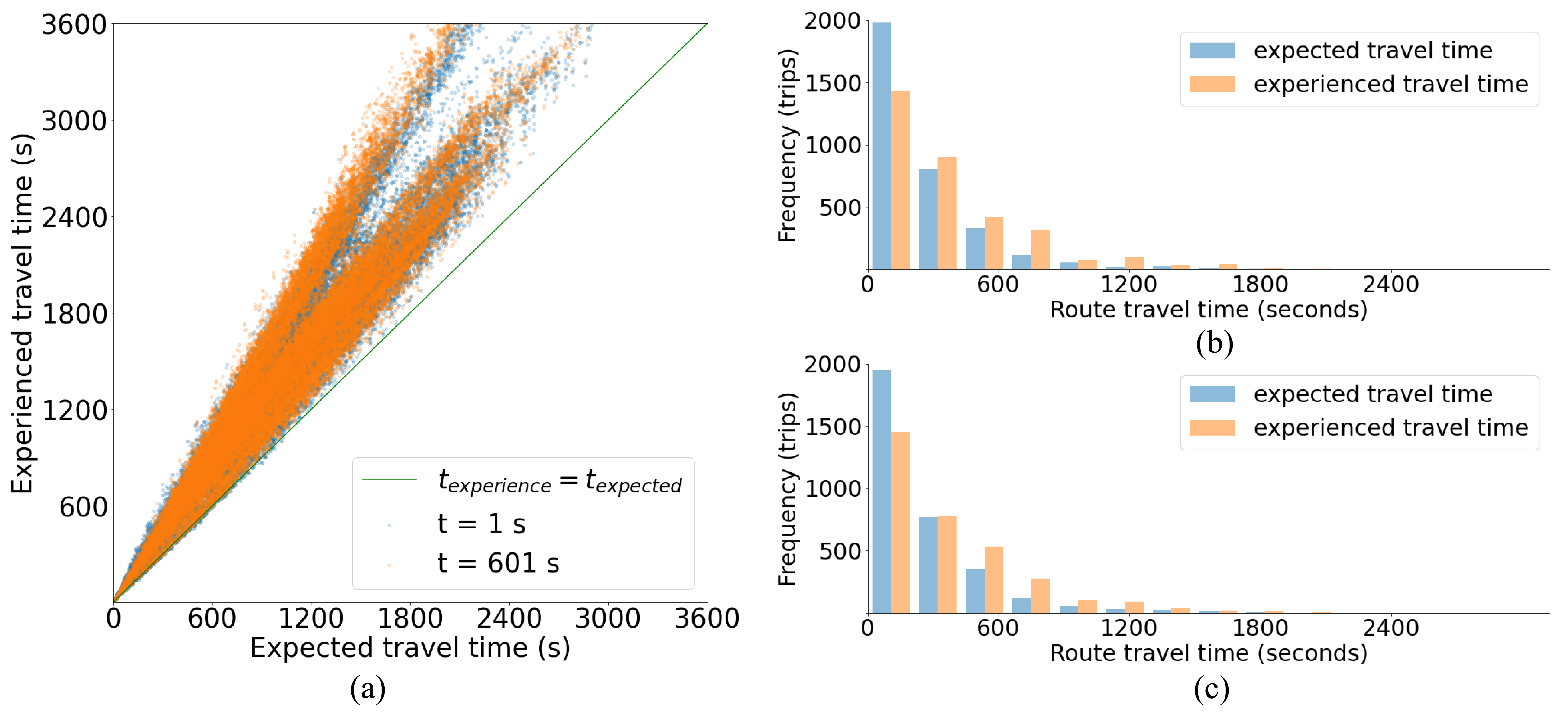}
    \caption{A comparison of the expected route travel times against experienced route travel times from the Sydney CBD footpath network: (a) Comparative scatter plot at $t$ = 1 sec and  $t$ = 601 sec, (b) distribution of route travel times between $t$ = 1 and $t$ = 300 sec, and (c) distribution of route travel times between $t$ = 601 sec and $t$ = 900 sec.}
    \label{fig:cbd_routett}
\end{figure}

The two are expected to have high correlation despite exhibiting slightly different patterns as shown in Figure \ref{fig:cbd_routett}. {\color{black} At the beginning of the simulation when the network is empty, experienced route travel times tend to be larger than the expected route travel times. However, 10 minutes into the simulation when the network becomes more congested, experienced route travel times tend to be even greater than the expected route travel times as shown in Figure \ref{fig:cbd_routett}(a). The distribution of experienced route travel times against expected route travel times are shown in Figure \ref{fig:cbd_routett}(b) for $t$ = 0 sec to $t$ = 600 sec and in Figure \ref{fig:cbd_routett}(c) for $t$ = 600 sec to $t$ = 1200 sec. The average experienced route travel time is greater than the expected route travel time. This results in a wider distribution of pedestrians across available paths in the network, preventing gridlocks to happen as a common limitation in many large-scale simulation-based DTA models \citep{mahmassani2013urban}. In Figure \ref{fig:cbd_routett}, we show two travel time distribution snapshots from $t$ = 1sec and $t$ = 601 sec to demonstrate how the expected and experienced route travel times change at different traffic states. The difference between the expected and experienced travel time at $t$ = 1 sec the $t$ = 601 sec does not appear to be significant. We believe that the 28,414 trips used as travel demand might not be too high to produce heavily congested traffic states.}

\raisebox{-0.45\height}

\FloatBarrier
\section{Conclusion} \label{sec:conclusion}
{\color{black} In this paper, we proposed a DPTA framework specifically for pedestrian networks with bidirectional links. The study formulated and presented the DUE problem for walking networks accounting for the microscopic characteristics of bidirectional pedestrian streams often observed in crowded sidewalks and crossings such as dynamic lane-formation and self-organization. To reproduce realistic formation and dissipation of shockwaves in the network, we presented an adaptation of the LTM consisting of a link model and a node model. The link model regulates how pedestrians traverse on any link, including how congestion due to bidirectional traffic affects the speed of pedestrians and how a capacity drop creates shockwaves in both directions. The link model also captures the bidirectional impact using a three-dimensional triangular bidirectional FD. The node model regulates how pedestrians go through any nodes following non-negativity, conservation of bidirectional flows and turning fractions constraints. 

We tested the developed DPTA modeling framework in a small grid network and a long corridor with strong bidirectional effects. Results from a small grid network successfully reproduced the pedestrians' behavior in avoiding routes with high travel time caused by the bidirectional effects and link closure. The simulated bidirectional long corridor successfully exhibited formation, propagation, and dissipation of shockwaves when link capacities change due to the bidirectional effects or introduction of a bottleneck. The paper also presented results from a large-scale sidewalk network model of Sydney CBD, demonstrating the applicability of the model to real-world networks.

Future research could further explore development of stochastic route choice models \citep{gentile2016solving}, adapting more efficient LTM algorithms \citep{raadsen2019continuous}, developing network-wide pedestrian control strategies \citep{molyneaux2021design} in the proposed DPTA context, and develop a route choice model under system optimum condition to improve overall safety and efficiency for large pedestrian crowd events.}

%% Adding acknolwedgment section
\section*{Acknowledgement}
This research was funded by the Australian Government through the Australian Research Council (project number DP220102382).

% Appendix here
% Options are (1) APPENDIX (with or without general title) or 
%             (2) APPENDICES (if it has more than one unrelated sections)
% Outcomment the appropriate case if necessary
%
% \begin{APPENDIX}{<Title of the Appendix>}
% \end{APPENDIX}
%
%   or 
%
% \begin{APPENDICES}
% \section{<Title of Section A>}
% \section{<Title of Section B>}
% etc
% \end{APPENDICES}

% References here (outcomment the appropriate case) 

% CASE 1: BiBTeX used to constantly update the references 
%   (while the paper is being written).
\bibliographystyle{informs2014trsc} % outcomment this and next line in Case 1
%\bibliography{<your bib file(s)>} % if more than one, comma separated
\bibliography{jabref}
% CASE 2: BiBTeX used to generate mypaper.bbl (to be further fine tuned)
%\input{mypaper.bbl} % outcomment this line in Case 2

\end{document}